%% file: 1997-030.tex
\documentclass[12pt,reqno]{amsart}
\usepackage{epsfig}

\let\epsfig\psfig

\usepackage{times,mathptm}

\title{On subdivision posets of cyclic polytopes}
\author{Paul H. Edelman}
\address{\hskip-\parindent Paul H. Edelman\\School of Mathematics\\University of Minnesota\\Minneapolis MN 55455}
\email{edelman@math.umn.edu}
\author{J\"org Rambau}\thanks{The second author was
  supported at MSRI in part by NSF grant \#DMS 9022140}
\address{\hskip-\parindent J\"org Rambau\\Konrad-Zuse-Zentrum f\"ur Informationstechnik\\
  Takustr.~7\\14195~Berlin\\
  Germany}
\email{rambau@zib.de}
\author{Victor Reiner}\thanks{The third author was supported by
  a Univ. of Minnesota McKnight-Land Grant Fellowship} 
\address{\hskip-\parindent Victor Reiner\\School of Mathematics\\University of Minnesota\\Minneapolis MN 55455}
\email{reiner@math.umn.edu}

\newcommand{\setof}[2]{\left\{\,{#1}\,:\,{#2}\:\right\}}
\newcommand{\romanenumi}%
{%
}
\newcommand{\red}{\mathrm{red}}
\newcommand{\green}{\mathrm{green}}
\newcommand{\Stwo}{\mathcal{S}_2}
\newcommand{\Sone}{\mathcal{S}_1}
\newcommand{\Sonetwo}{\mathcal{S}}
\newcommand{\proper}[1]{\overline{#1}}
\newcommand{\sm}{\backslash}
\newcommand{\lk}{\operatorname{lk}}
\newcommand{\st}{\operatorname{st}}
\newcommand{\as}{\operatorname{del}}
\newcommand{\conv}{\operatorname{conv}}
\newcommand{\step}[1]{\par \emph{Step #1:}}
\newcommand{\case}[1]{\par \textsc{Case #1:}}
\newtheorem{Theorem}{Theorem}[section]

\newtheorem{Conjecture}[Theorem]{Conjecture}
\newtheorem{Lemma}[Theorem]{Lemma}
\newtheorem{Sublemma}[Theorem]{Sublemma}
\newtheorem{Proposition}[Theorem]{Proposition}

\newcommand{\oo}{\operatorname{Baues}(C(n,d))}

\newcommand{\susp}{\operatorname{susp}}
\newcommand{\Int}{\operatorname{Int}}
\newcommand{\intbar}{\proper{\operatorname{Int}}}
\newcommand{\atomic}{\mathit{atomic}}
\newcommand{\coatomic}{\mathit{coatomic}}
\newcommand{\submersion}{\operatorname{sub}_{\lceil\frac{d}{2}\rceil}}

\begin{document}

\begin{abstract}
  There are two related poset structures, the higher Stasheff-Tamari
  orders, on the set of all triangulations
  of the cyclic $d$ polytope with $n$ vertices.
  In this paper it is shown that both of them
  have the homotopy type of a sphere of dimension $n-d-3$. 

  Moreover, we resolve
  positively a new special case of the
  \emph{Generalized Baues Problem}:
  The Baues poset of all polytopal decompositions of a cyclic
  polytope of dimension $d \leq 3$ has the homotopy type of a sphere
  of dimension $n-d-2$.
\end{abstract}

\maketitle

%%%%%%%%%%%%%%%%%%%%%%%%%%%%%%%%%%%%%%%%%%%%%%%%%%%%%%%%%%%%%%%%%%%%%%%%%%%%%%%
%% Introduction:                                                             %%
%% 07/03/97 pe/jr/vr                                                         %%
%%%%%%%%%%%%%%%%%%%%%%%%%%%%%%%%%%%%%%%%%%%%%%%%%%%%%%%%%%%%%%%%%%%%%%%%%%%%%%%
\section{Introduction}

This paper continues the investigation of certain posets of triangulations
of cyclic polytopes, the \emph{higher Stasheff Tamari posets}, initiated in
\cite{EdelmanReiner1996} and continued in~\cite{Rambau1996}. 

The first higher Stasheff Tamari poset is the poset $\Sone(n,d)$ of all
triangulations of the cyclic $d$-polytope with $n$ vertices $C(n,d)$, partially
ordered by increasing bistellar operations;
the second higher Stasheff Tamari poset is the poset $\Stwo(n,d)$ of all
triangulations of $C(n,d)$, partially ordered by the height of their
characteristic sections in $C(n,d+1)$
(see \cite{EdelmanReiner1996}, \cite{Rambau1996}).

Our first main result is the following.

\begin{Theorem}
  \label{thm:introduction:hstmain}
  \begin{enumerate}\romanenumi
  \item
  For all $n > d + 1$ the proper part $\proper{\Sone(n,d)}$
  of $\Sone(n,d)$ is homotopy equivalent to a sphere of dimension $n-d-3$.
  \item
  For all $n > d + 1$ the proper part $\proper{\Stwo(n,d)}$
  of $\Stwo(n,d)$ is homotopy equivalent to a sphere of dimension $n-d-3$.
  \end{enumerate}
\end{Theorem}

In \cite{EdelmanReiner1996}, it was proved for $d \leq 3$ that the
poset structures $\Sone(n,d)$ and $\Stwo(n,d)$ coincide. It was also
shown that the poset $\Stwo(n,d)$ is a lattice
for $d \leq 3$. If $d=2$ this is the well-known \emph{Tamari lattice}
on triangulations of a convex $n$-gon.  We will use this lattice structure
to resolve in the affirmative a special case of the
\emph{Generalized Baues Problem} of Billera, Kapranov, and Sturmfels (see
\cite{BilleraKapranovSturmfels1994},
\cite{Rambau1996a},~\cite{RambauZiegler1995}).

\begin{Theorem}
  \label{thm:introduction:baues3main}
  For cyclic polytopes $C(n,d)$ of dimension $d \leq 3$,
  the refinement ordering on the set of polytopal subdivisions  
  gives a poset which is homotopy equivalent to a $(n-d-2)$-sphere.
\end{Theorem}

We will prove Theorem~\ref{thm:introduction:hstmain} by induction
on $n-d$, showing that the poset $\proper{\Sone(n,d)}$ (resp.\
$\proper{\Stwo(n,d)}$) is homotopy equivalent to the suspension of
$\proper{\Sone(n-1,d)}$ (resp.\ $\proper{\Stwo(n-1,d)})$.

The proof of Theorem~\ref{thm:introduction:baues3main} is via a reduction
to the poset $\proper{\Stwo(n,d)}$, by showing that
the poset of polytopal subdivisions of $C(n,d)$ is homotopy equivalent to the
suspension of $\proper{\Stwo(n,d)}$. We will make use of a lemma
(Lemma~\ref{Webb-lemma})
about the homotopy type of non-contractible intervals in a poset
which we think is of interest in its own right.

This paper is structured as follows: in Section~\ref{sec:basics} we
recall some notation and basic facts about simplicial complexes,
posets, and cyclic polytopes. In Section~\ref{sec:homotopy} we
prove Theorem~\ref{thm:introduction:hstmain}. Sections~\ref{sec:lemmasS1}
and~\ref{sec:lemmasS2} provide the necessary details.
%thereby
%revealing a lot of properties of triangulations of cyclic polytopes
%that we consider interesting in their own right.
In Section~\ref{sec:baues} we
prove Theorem~\ref{thm:introduction:baues3main}, the special case of the
Generalized Baues Problem.
Section~\ref{sec:problems} discusses some of the remaining
open problems in the area of triangulations of cyclic polytopes.

%%%%%%%%%%%%%%%%%%%%%%%%%%%%%%%%%%%%%%%%%%%%%%%%%%%%%%%%%%%%%%%%%%%%%%%%%%%%%%%
%% Notation and Basic facts:                                                 %%
%% 07/03/97 pe/jr/vr                                                         %%
%%%%%%%%%%%%%%%%%%%%%%%%%%%%%%%%%%%%%%%%%%%%%%%%%%%%%%%%%%%%%%%%%%%%%%%%%%%%%%%
\section{Notation and Basic Facts}
\label{sec:basics}

In this section we will introduce our notation and discuss some basic
facts that have appeared previously.

Let $[n]:=\{1,2,\ldots,n\}$.
We regard the \emph{cyclic $d$-polytope with $n$ vertices} as the convex hull of
points on the moment curve 
\begin{displaymath}
  C(n,d) = \conv \{ (i, i^2, \dots , i^d) \in \mathbb{R}^d : i \in [n] \}
\end{displaymath}
Since we are dealing with the combinatorial structure of all
triangulations of
cyclic polytopes we may choose  these
special coordinates without  any loss of generality.  We will often
refer to the $i^{\mathit{th}}$ vertex $(i,i^2,\ldots,i^d)$ of $C(n,d)$
as simply $i$.

The \emph{canonical projection} $p = p_{n,d}$ from $C(n, d+1)$
onto $C(n, d)$ is
given by deletion of the $x_{d+1}$-coordinate.
Facets of $C(n,d)$ that can be seen from a point in $\mathbb{R}^{d+1}$
with a very large (negative) $x_{d+1}$-coordinate
are called \emph{upper (lower) facets}.
 
Two simplices are said to be \emph{admissible} if they 
intersect in a common (possibly empty) face of each.
A \emph{triangulation} of a polytope~$P$ is a set of simplices with
vertices in the vertex set of~$P$ such that
\begin{itemize}
\item the union of the simplices equals~$P$,
\item every face of a simplex in the triangulation is itself
  in the triangulation, and
\item any two simplices are admissible.
\end{itemize}

Triangulations are often identified with their sets of inclusion-maximal
faces. Simplices are usually identified with their vertex sets.

To test intersections
of simplices $S_1$ and $S_2$
we will use the concept of \emph{zig-zag-paths} based on
the alternating oriented matroid property of cyclic polytopes
(see~\cite{Rambau1996}).
We construct a table with $n$ columns, corresponding to the
labels $1, \dots , n$, and two rows, corresponding to the simplices $S_1$
and~$S_2$. In row $i$, column $j$, there is a star $\ast$ if and only if
$j \in S_i$.
An \emph{$(S_1,S_2)$-zig-zag-path of length~$k$} is a set of $k$ stars in the
columns $s_1 < \dots < s_k$
such that $s_1, s_3, s_5, \dots$ are in $S_1$ and $s_2, s_4, s_6, \dots$
are in $S_2$, or vice versa.
The simplices $S_1$ and $S_2$ are admissible in dimension~$d$ if and only if
there is no $(S_1,S_2)$-zig-zag-path of length~$d+2$
(see Figure~\ref{fig:notation:zigzag}).

\begin{figure}[t]
  \begin{center}
    \leavevmode
    \input{bauclic_zigzagexample.pstex_t}
    \caption{Zig-zag-paths: (a) $S_1$ and $S_2$ are non-admissible
      in dimensions $5$ or less, hence $S_1$ and $S_2$ cannot be in a
      triangulation of, e.~g., $C(9,4)$ at the same time;
      (b) $S_1$ and $S_2$ are admissible in dimensions $3$ or greater,
      therefore $S_1$ and $S_2$ may be in a triangulation of $C(9,4)$.}
    \label{fig:notation:zigzag}
  \end{center}
\end{figure}

Any subset $V \subseteq [n]$ gives rise to a cyclic subpolytope $C(V,d)$,
namely the convex hull of the subset $V$.
For a cyclic subpolytope $C(V,d)$ and
a $d$-subset $F$ of $V$ we call a label $i \in V \sm F$
an \emph{even (odd) gap} of $F$ (in $V$) if the number of labels $j \in V$
with $j > i$ is even (odd).
Then we know that the set of \emph{lower (upper) facets}
of $C(V,d)$ is the set of all $F \in \tbinom{V}{d}$
containing only even (odd)
gaps~\cite{Rambau1996}. This applies, in particular, to simplices so
that we can talk about the upper and lower facets of a $d$-simplex
in $C(n,d)$.
For a visualization, we use the same table as for the zig-zag-paths and
fill an $e$ (resp.\ $o$) for an even (resp.\ odd) gap into the corresponding
field (see Figure~\ref{fig:notation:gap}).

\begin{figure}[t]
  \begin{center}
    \leavevmode
    \input{bauclic_gapexample.pstex_t}
    \caption{Gaps: (a) $S_1$ is a lower facet of $C(9,4)$,
      $S_2$ is an upper facet of $C(9,4)$;
      (b) neither $S_1$ nor $S_2$ are facets of $C(9,4)$, but
      $S_1$ is a lower facet and $S_2$ an upper facet of the
      cyclic subpolytope $C(\{ 2,3,6,7,8,9 \},4)$.
      Therefore, $S_1$ is lower than $S_2$ in $C(9,5)$.}
    \label{fig:notation:gap}
  \end{center}
\end{figure}

Let $T$ be a triangulation of $C(n,d)$ and $\Tilde{S}$ be a $(d+1)$-simplex
in $C(n,d+1)$ all of whose lower facets lie in $T$.
An \emph{increasing bistellar operation} or \emph{increasing flip} 
in $T$ at $\Tilde{S}$ is
an operation that replaces in $T$ the lower facets of $\Tilde{S}$
by the upper facets of $\Tilde{S}$.  The result, it is clear, is
a new triangulation of $C(n,d)$.
The transitive closure of this operation defines the
\emph{first higher Stasheff-Tamari poset} $\Sone(n,d)$.
We write $T <_1 T'$ to indicate that $T$ is less than $T'$ in $\Sone(n,d)$.

The \emph{characteristic section} of a triangulation $T$ of $C(n,d)$ is
the unique piecewise linear map (with respect to the simplicial complex~$T$)
from $C(n,d)$ to $C(n,d+1)$ that is inverted by the canonical projection~$p$
and has the property that it sends the $i^{\mathit{th}}$ vertex of $C(n,d)$ to
the $i^{\mathit{th}}$ vertex of $C(n,d+1)$.
We identify a triangulation $T$ with its characteristic section
$T: C(n,d) \to C(n,d+1)$ and with its image $T(C(n,d))$ in
$C(n,d+1)$. The \emph{second higher Stasheff-Tamari poset} $\Stwo(n,d)$ is the
set of all triangulations of $C(n,d)$ partially ordered by the
height of characteristic sections. That is, $T \le_2 T'$ if and only
if $T(x)_{d+1} \le T'(x)_{d+1}$ for all $x \in C(n,d)$, where here
$v_{d+1}$ denotes the $(d+1)^{\mathit{st}}$ coordinate of the
vector $v$ in $\mathbb{R}^{d+1}$.
We then say that $T$ is weakly lower than $T'$. If $T(x)_{d+1} \le T'(x)_{d+1}$
holds for all $x$ in the (geometric) intersection of a simplex
$S \in T$ and a simplex
$S' \in T'$ we say that $S$ is weakly lower than $S'$.
We write $T <_2 T'$ to denote that $T$ is less than $T'$ in $\Stwo(n,d)$.

The unique minimal element in $\Sone(n,d)$ respectively $\Stwo(n,d)$
(which is the set of lower facets of $C(n,d+1)$) is denoted
by $\Hat{0}_{n,d}$. Similarly, the unique maximal element 
(which is the set of upper facets of $C(n,d+1)$) is denoted by $\Hat{1}_{n,d}$.
The $d$-simplices in $C(n,d)$ are partially ordered by
the following relation: $S \prec S'$ if and only if $S \cap S'$ is
a lower facet of $S'$ and an upper facet of~$S$ (see~\cite{Rambau1996}).

We will make use of some standard constructions on simplicial
complexes.  Let $\Delta$ be a simplicial complex on the ground set
$X$.  That is, $\Delta$ is a collection of subsets  of $X$ that
is closed under containment.  If $S \subseteq X$ define the \emph{
link of $S$ in $\Delta$} to be the complex
\begin{displaymath}
  \lk_{\Delta}(S):=\setof{R \sm S}{R \in \Delta, S \subseteq R};
\end{displaymath}
the \emph{star of $S$ in $\Delta$} is the complex
\begin{displaymath}
  \st_{\Delta}(S):=\setof{R \in \Delta}{S \subseteq R};
\end{displaymath}
and the \emph{deletion of $S$ in $\Delta$} is the complex
\begin{displaymath}
  \as_{\Delta}(S):=\setof{R \in \Delta}{S \not\subseteq R}.
\end{displaymath}
If there is another complex $\Delta'$ on a ground set $Y$ disjoint from
$X$ we will define the \emph{combinatorial join of $\Delta$ and $\Delta'$}
to be the complex on the ground set $X \cup Y$
\begin{displaymath}
  \Delta * \Delta' := \setof{S \cup S'}{S \in \Delta , S'\in \Delta'}.
\end{displaymath}

If $T, T'$ are the sets of inclusion maximal faces of $\Delta , \Delta'$
then the above formulas yield the sets of inclusion maximal faces of
the link, the star, the deletion, and the join, respectively.

Given an $i$-simplex $\sigma$ spanned by some $(i+1)$-subset
(also denoted $\sigma$) of vertices of $C(n,d)$, there is also
a unique linear section $\sigma: \sigma \rightarrow C(n,d+1)$
of $p$ having the property that it sends each vertex $i$ of $\sigma$ to
the vertex labelled $i$ of $C(n,d+1)$.
Say that $\sigma$ \emph{submerged} by the triangulation $T$ of $C(n,d)$ 
if
\begin{displaymath}
  \sigma(x)_{d+1} \leq T(x)_{d+1}
\end{displaymath}
for every point $x$ in $\sigma$.  For a triangulation $T$ of $C(n,d)$
let its $i^{th}$ \emph{submersion set} $\operatorname{sub}_i(T)$ be
the set of $i$-simplices submerged by $T$.

When we refer to the topology or homotopy type of a poset $P$,
we will always mean the topology of the \emph{geometric realization}
of its \emph{order complex}, i.~e.,
$\vert \Delta(P) \vert$~\cite[\S9]{Bjoerner1995}.
If $P$ is a poset with bottom and top elements $\hat{0}, \hat{1}$,
then its \emph{proper part} $\proper{P}$ is simply the subposet
$P \sm \{\hat{0},\hat{1}\}$.

We recall the following facts from \cite{EdelmanReiner1996}
and~\cite{Rambau1996}
which will be crucial for our main results:

\begin{Theorem}\cite[Theorem 1.1]{Rambau1996}
  \label{thm:basics:bounded}
  The first higher Stasheff Tamari poset $\Sone(n,d)$ is bounded.
\end{Theorem}

\begin{Theorem}\cite[Theorem 4.2(iii), Proposition 5.14(iii)]{Rambau1996}
  \label{thm:basics:f}
  The following map is well-defined and order-preserving:
  \begin{displaymath}
  f:
  \left\{
  \begin{array}{rcl}
    \Sone(n,d) & \to & \Sone(n-1,d),\\
    T & \mapsto & \as_T(n) \cup \left(\as_{\lk_T(n)}(n-1) * \{n-1\}\right).
  \end{array}
  \right.    
  \end{displaymath}
\end{Theorem}

\begin{Proposition}
\label{submersion}\cite[Proposition 2.15 ]{EdelmanReiner1996}
For any two triangulations $T_1,T_2$ of $C(n,d)$, we have
$T_1 \leq T_2$ in $\Stwo(n,d)$ if and only if 
\begin{displaymath}
  \submersion(T_1) \subseteq \submersion(T_2).
\end{displaymath}
\end{Proposition}

\begin{Proposition}
\label{low-dimensional-submersion}
\cite[Propositions 3.2, 4.1 ]{EdelmanReiner1996}
Membership in $\lceil\frac{d}{2}\rceil$-
submersion sets for $d=2, 3$ has the following
characterization.

  For $T$ a triangulation of $C(n,2)$ and $e=\{i,j\}$ an
edge inside $C(n,2)$, we have that $e \in sub_1(T)$ if and only
if there does not exist an edge $e'=\{k,l\}$ of $T$ with $k<i<l<j$.

  For $T$ a triangulation of $C(n,3)$ and $t=\{i,j,k\}$ a
triangle inside $C(n,3)$, we have that $t \in sub_2(T)$ if and only
if there does not exist an edge $\{x,y\}$ of $T$ with $i<x<j<y<z$. 
\end{Proposition}

\begin{Theorem}\cite[Theorems 3.6, 4.9]{EdelmanReiner1996}
\label{latticeness}
For $d \leq 3$, the higher Stasheff-Tamari poset $\Stwo(n,d)$ is
a lattice, i.~e., any subset of its elements has a meet (greatest
lower bound) and a join (least upper bound).
\end{Theorem}

\begin{Theorem}\cite[Theorems 3.9, 4.11]{EdelmanReiner1996}
  \label{Tamari-sphericity}
  For $d \leq 3$, the proper part $\proper{\Stwo(n,d)}$ of
  the higher Stasheff-Tamari poset has the homotopy type of an
  $(n-d-3)$-sphere.
\end{Theorem}

%%%%%%%%%%%%%%%%%%%%%%%%%%%%%%%%%%%%%%%%%%%%%%%%%%%%%%%%%%%%%%%%%%%%%%%%%%%%%%%
%% The homotopy types of $\Sone(n,d)$ and $\Stwo(n,d)$:                      %%
%% 07/03/97 pe/jr/vr                                                         %%
%%%%%%%%%%%%%%%%%%%%%%%%%%%%%%%%%%%%%%%%%%%%%%%%%%%%%%%%%%%%%%%%%%%%%%%%%%%%%%%

\section{The Homotopy Types of $\Sone(n,d)$ and $\Stwo(n,d)$}
\label{sec:homotopy}

In this section, Theorem~\ref{thm:introduction:hstmain} will be proven
by induction on $n-d$, using the Suspension
Lemma~\ref{thm:introduction:suspensionlemma} below to
show that the proper part of $\Sonetwo(n,d)$ is homotopy equivalent to
the suspension of the proper part of $\Sonetwo(n-1,d)$,
where $\Sonetwo(n,d)$ can be either $\Sone(n,d)$ or $\Stwo(n,d)$.
(A more detailed proof of the Suspension Lemma can be found
in~\cite{Rambau1997}.)

%%\begin{Lemma}[Suspension Lemma]
%%\label{thm:introduction:suspensionlemma}
%%Let $P,Q$ be bounded posets with $\Hat{0}_Q \neq \Hat{1}_Q$. Moreover, let
%%\begin{displaymath}
%%  f: P \to Q \quad \text{and} \quad i,j: Q \to P
%%\end{displaymath}
%%be order-preserving with the following properties:
%%\begin{enumerate}\romanenumi
%%\item \label{itm:susp:Pgreenred}
%%  The elements of $P$ are the disjoint union of $\green(P)$ 
%%  (the \emph{green} elements of $P$) and
%%  $\red(P)$ (the \emph{red} elements of $P$)
%%  such that $\green(P)$ is an order ideal in~$P$ (or, equivalently,
%%  $\red(P)$ is an order filter in~$P$).
%%\item \label{itm:susp:composition}
%%  The maps $f \circ i$ and $f \circ j$ are the identity on~$Q$.
%%\item \label{itm:susp:i00j11}
%%  The map $i$ maps $\Hat{0}_Q$ to $\Hat{0}_P$, and
%%  $j$ maps $\Hat{1}_Q$ to $\Hat{1}_P$.
%%\item \label{itm:susp:igreenjred}
%%  The image of $i$ is green, and the image of $j$ is red.
%%\item \label{itm:susp:interval}
%%  For every $q \in Q$ the preimage $f^{-1}(q)$ is the interval
%%  $[i(q),j(q)]$.
%%\item \label{itm:susp:red0green1}
%%  The interval $[i(\Hat{0}_Q),j(\Hat{0}_Q)]$ is red except for $\Hat{0}_P$,
%%  and the interval $[i(\Hat{1}_Q),j(\Hat{1}_Q)]$ is green except
%%  for $\Hat{1}_P$.
%%\end{enumerate}
%%Then the proper part $\proper{P}$ of $P$ is homotopy equivalent to the
%%suspension of the proper part $\proper{Q}$ of~$Q$.
%%\end{Lemma}

\begin{Lemma}[Suspension Lemma]
\label{thm:introduction:suspensionlemma}
Let $P,Q$ be bounded posets with $\Hat{0}_Q \neq \Hat{1}_Q$. Assume
there exist a dissection of $P$ into green elements $\green(P)$ and
red elements $\red(P)$, as well as order-preserving maps
\begin{displaymath}
  f: P \to Q \quad \text{and} \quad i,j: Q \to P
\end{displaymath}
with the following properties:
\begin{enumerate}\romanenumi
\item \label{itm:susp:Pgreenred}
  The green elements form an order ideal in~$P$.
%%  (or, equivalently, the red elements form an order filter in~$P$).
\item \label{itm:susp:composition}
  The maps $f \circ i$ and $f \circ j$ are the identity on~$Q$.
\item \label{itm:susp:igreenjred}
  The image of $i$ is green, the image of $j$ is red.
\item \label{itm:susp:interval}
  For every $p \in P$ we have $(i \circ f)(p) \le p \le (j \circ f)(p)$.
\item \label{itm:susp:red0green1}
  The fiber $f^{-1}(\Hat{0}_Q)$ is red except for $\Hat{0}_P$, the fiber
  $f^{-1}(\Hat{1}_Q)$ is green except for~$\Hat{1}_P$.
\end{enumerate}
Then the proper part $\proper{P}$ of $P$ is homotopy equivalent to the
suspension of the proper part $\proper{Q}$ of~$Q$.
\end{Lemma}

\begin{proof}[Sketch of proof.]
  Define
\begin{align}
  g: 
  &\left\{
  \begin{array}{rcl}
    \proper{P} & \to     & \proper{Q \times \{ \Hat{0}, \Hat{1} \}},\\
    p       & \mapsto &
    \begin{cases}
      (f(p), \Hat{0}) & \text{if $p$ is green},\\
      (f(p), \Hat{1}) & \text{if $p$ is red};
    \end{cases}
  \end{array}
  \right.\\
\intertext{and}
  h: 
  &\left\{
  \begin{array}{rcl}
    \proper{Q \times \{ \Hat{0}, \Hat{1} \}} & \to & \proper{P},\\
    (q, \Hat{0}) & \mapsto & i(q),\\
    (q, \Hat{1}) & \mapsto & j(q).
  \end{array}
  \right.
\end{align}

The assumptions guarantee that the above maps are well-defined and
order-preserving.
Observe that $g \circ h$ is the identity map
on $\proper{Q \times \{ \Hat{0}, \Hat{1} \}}$ and that
$\proper{Q \times \{ \Hat{0}, \Hat{1} \}}$ is homeomorphic to the suspension
of~$\proper{Q}$. It is easy to show
that both $h \circ g$ and the identity map on~$\proper{P}$ are carried by the
following contractible carrier on the order complex $\Delta (\proper{P})$
of~$\proper{P}$.

%%\begin{displaymath}
%%%  C: 
%%%  \left\{
%%   \begin{array}{rcl}
%%%    \Delta(\proper{P}) & \to & \Delta(\proper{P}),\\
%%    \sigma & \mapsto &
%%    \begin{cases}
%%      \proper{P}_{\ge (i \circ f)(\min \sigma)} &
%%      \text{if $f(\min \sigma) > \Hat{0}_Q,
%%        f(\max \sigma) = \Hat{1}_Q$},\\
%%      \proper{P}_{\le (j \circ f)(\max \sigma)} &
%%      \text{if $f(\max \sigma) < \Hat{1}_Q,
%%        f(\min \sigma) = \Hat{0}_Q$},\\
%%      \proper{P}_{\ge (i \circ f)(\min \sigma)}\\
%%      {} \cap
%%      \proper{P}_{\le (j \circ f)(\max \sigma)} &
%%      \text{if $f(\min \sigma) > \Hat{0}_Q,
%%        f(\max \sigma) < \Hat{1}_Q$}.
%%    \end{cases}
%%  \end{array}
%%%  \right.
%%\end{displaymath}

\begin{displaymath}
  C: 
  \left\{
  \begin{array}{rcl}
    \Delta(\proper{P}) & \to & 2^{\Delta(\proper{P})},\\
    \sigma & \mapsto &
      \Delta
      \bigl(
      P_{\ge (i \circ f)(\min \sigma)}
      \cap
      P_{\le (j \circ f)(\max \sigma)}
      \cap
      \proper{P}
      \bigr).
  \end{array}
  \right.
\end{displaymath}

Thus, by the Carrier Lemma~\cite{Bjoerner1995}, the map $h \circ g$ is
homotopic to the identity on $\proper{P}$, and $g$ and $h$ are homotopy
inverses to each other.
\end{proof}

We now prove that the assumptions of
the Suspension Lemma are satisfied by the following set of data.

\begin{align*}
  P &= \Sonetwo(n,d),\\
  Q &= \Sonetwo(n-1,d),\\
  \green(\Sonetwo(n,d)) &= 
  \begin{cases}
    \setof{T \in \Sonetwo(n,d)}{
    \{ n-d, \dots , n \} \notin T} & \text{for $d$ even},\\
    \setof{T \in \Sonetwo(n,d)}{
    \{ n-d, \dots , n \} \in T}    & \text{for $d$ odd};
  \end{cases}\\
  \red(\Sonetwo(n,d)) &= 
  \begin{cases}
    \setof{T \in \Sonetwo(n,d)}{
    \{ n-d, \dots , n \} \in T}    & \text{for $d$ even},\\
    \setof{T \in \Sonetwo(n,d)}{
    \{ n-d, \dots , n \} \notin T} & \text{for $d$ odd};
  \end{cases}\\
  f:
  &\left\{
  \begin{array}{rcl}
    \Sonetwo(n,d) & \to & \Sonetwo(n-1,d),\\
    T & \mapsto & \as_T(n) \cup \as_{\lk_T(n)}(n-1) * \{n-1\};
  \end{array}
  \right.\\
  i:
  &\left\{
  \begin{array}{rcl}
    \Sonetwo(n-1,d) & \to & \Sonetwo(n,d),\\
    T & \mapsto & T \cup \st_{\Hat{0}_{n,d}}(n);
  \end{array}
  \right.\\
  j:
  &\left\{
  \begin{array}{rcl}
    \Sonetwo(n-1,d) & \to & \Sonetwo(n,d),\\
    T & \mapsto & \as_T(n-1)\\
    & & {} \cup \lk_T(n-1) * \{n\}\\
    & & {} \cup \st_{\Hat{1}_{n,d}}(\{n-1,n\}).
  \end{array}
  \right.
\end{align*}

Theorem~\ref{thm:basics:bounded} shows that $\Sonetwo(n,d)$ is
bounded. Moreover, by Theorem~\ref{thm:basics:f} we know that $f(T)$ is a
triangulation of $C(n-1,d)$ for all triangulations $T$ of
$C(n,d)$.  The geometric description of $f$ is as follows: 
starting with the triangulation $T$ of $C(n,d)$, if one slides the
vertex $n$ along the moment curve until it coincides with the vertex
$n-1$, then certain $d$-simplices of $T$ will degenerate.  Removing these
degenerate simplices and renaming all occurrences of $n$ by $n-1$ yields
the triangulation $f(T)$.
 
The constructions of $i$ and $j$ can be described
geometrically as follows: The cyclic polytope $C(n-1,d)$ can be
embedded into the cyclic polytope $C(n,d)$ in many different ways.
For example, there is an embedding that sends vertex $k$ of $C(n-1,d)$ to
vertex $k$ in $C(n,d)$ for all $1 \le k \le n-1$.
There is another embedding which sends vertex $k$ to vertex $k$
for all $1 \le k < n-1$ and vertex $n-1$ of $C(n-1,d)$ to vertex $n$ of
$C(n,d)$.

The map $i$ uses the first embedding of $C(n-1,d)$ into $C(n,d)$ to
embed a triangulation $T$ of $C(n-1,d)$ into $C(n,d)$. This leads to
a partial triangulation of $C(n,d)$. Since the
``new'' vertex $n$ in $C(n,d)$ ``sees'' a convex polytope from
outside, the only
possibility to complete that partial triangulation 
is to join every facet of $T$ that is ``visible'' by $n$
to~$n$. It is an easy calculation using Gale's Evenness Criterion
\cite[Theorem~0.7]{Ziegler1994}
that the given formula for $i$ describes exactly that.

The map $j$ uses the second embedding of $C(n-1,d)$ into $C(n,d)$ for
embedding a triangulation $T$ of $C(n-1,d)$ into $C(n,d)$.
Again, the ``new'' vertex $n-1$ ``sees'' certain facets of a cyclic
polytope with $n-1$ vertices. Given a triangulation of $C(n-1,d)$ that
is embedded into $C(n,d)$ in this fashion, the only way to complete it to
a triangulation of $C(n,d)$ is to join $n-1$ with the visible
facets of the embedded $C(n-1,d)$.
Gale's Evenness Criterion again allows us to obtain the formula
for~$j$.
This proves that $i$ and $j$ are well-defined.

In the following we
outline the proof of Theorem~\ref{thm:introduction:hstmain} by verifying
the assumptions of the Suspension Lemma. Whenever the details
are more involved we give a reference to a Lemma in
Section~\ref{sec:lemmasS1} or \ref{sec:lemmasS2}, respectively.

If $n > d+2$ then $\Hat{0}_{n-1,d} \neq \Hat{1}_{n-1,d}$.
From Lemma \ref{thm:lemmasS1:maps} (resp.~\ref{thm:lemmasS2:maps}) we get
that all maps are order-preserving.
From Lemma \ref{thm:lemmasS1:specialsimplex}
(resp.~\ref{thm:lemmasS2:specialsimplex}) we know that no green element
can be above a red one.
By construction, $f \circ i$ and  $f \circ j$ are both
the identity on $\Sonetwo(n,d)$.
Since whether or not $\{ n-d, \dots , n \}$ is contained in $i(T)$ (resp.\
$j(T)$) does not depend on $T$, it can easily be seen that the image
of $i$ is green and that the image of $j$ is red.
From Lemma \ref{thm:lemmasS1:composition} (resp.~\ref{thm:lemmasS2:composition})
it follows that the preimages of any $T \in \Sonetwo(n,d)$ under $f$
are bounded by $i(T)$ and $j(T)$.
Finally, Lemma \ref{thm:lemmasS1:fibers01}
(resp.~\ref{thm:lemmasS2:fibers01}) imply that $\Hat{0}_{n,d}$ is the only green
element in $f^{-1}(\Hat{0}_{n-1,d})$ and that
$\Hat{0}_{n,d}$ is the only red element in $f^{-1}(\Hat{1}_{n-1,d})$.

The proof of Theorem~\ref{thm:introduction:hstmain} then follows from
the well-known fact that $C(d+2,d)$ has exactly two triangulations (i.~e.,
its proper part is the empty set which is a $(-1)$-sphere)
and induction on the codimension $n-d$ using
the Suspension Lemma~\ref{thm:introduction:suspensionlemma}.

%%%%%%%%%%%%%%%%%%%%%%%%%%%%%%%%%%%%%%%%%%%%%%%%%%%%%%%%%%%%%%%%%%%%%%%%%%%%%%%
%% Lemmas for $\Sone(n,d)$:                                                  %%
%% 07/03/97 pe/jr/vr                                                         %%
%%%%%%%%%%%%%%%%%%%%%%%%%%%%%%%%%%%%%%%%%%%%%%%%%%%%%%%%%%%%%%%%%%%%%%%%%%%%%%%

\section{Lemmas on $\Sone(n,d)$}
\label{sec:lemmasS1}

We first formulate a lemma that we are using to establish the comparability
of elements in $\Sone(n,d)$.

\setcounter{Theorem}{-1}

\begin{Lemma}
  \label{thm:lemmasS1:S1}
  Let $T$ and $T'$ be triangulations of $C(n,d)$.
  $T$ is less than or equal to $T'$ in $\Sone(n,d)$ if and only if
  there is a triangulation of the region between the characteristic sections
  of $T$ and $T'$ in $C(n,d+1)$.

  In other words, $T \le_1 T'$ if and only if
  there is a
  set $\Tilde{T}$ of $(d+1)$-simplices such that
  the following hold:
  \begin{enumerate}\romanenumi
  \item 
    Every pair of $(d+1)$-simplices in $T$ are admissible.
  \item
    For every lower facet $S$ of a $(d+1)$-simplex
    in $\Tilde{T}$ either there is another $(d+1)$-simplex in $\Tilde{T}$
    containing $S$, or $S$ is in~$T$. 
  \item
    For every upper facet $S$ of a $(d+1)$-simplex
    in $\Tilde{T}$ either there is another $(d+1)$-simplex in $\Tilde{T}$
    containing $S$, or $S$ is in~$T'$.
  \item Every $d$-simplex in $T \sm T'$ is a lower facet of some
    $(d+1)$-simplex in~$\Tilde{T}$.
  \item Every $d$-simplex in $T' \sm T$ is an upper facet of some
    $(d+1)$-simplex in~$\Tilde{T}$.
  \item Every $d$-simplex in $T \sm T' \cup T' \sm T$ is a facet of at most one
    $(d+1)$-simplex in $\Tilde{T}$.
  \end{enumerate}

  If the above assumptions are met we say
  ``$\Tilde{T}$ connects $T$ and~$T'$.''
\end{Lemma}

\begin{proof}
  Given a set of $(d+1)$-simplices $\Tilde{T}$ as in the assumption
  we get a sequence   of increasing flips from $T$ to $T'$ by
  sorting the simplices of $\Tilde{T}$ by any linear
  extension of ``$\prec$,'' as was
  shown in~\cite{Rambau1996}. On the other hand, every set of
  $(d+1)$-simplices corresponding to a sequence of
  increasing flips from $T$ to $T'$ has the properties listed above.
\end{proof}

We now prove a sequence of lemmas that allows us to apply the
Suspension Lemma in the case of $\Sone(n,d)$. Throughout this section it is
always assumed that $n > d+2$.

\begin{Lemma}
  \label{thm:lemmasS1:maps}
  The following maps are order-preserving.
  \begin{align*}
    f:
    &\left\{
      \begin{array}{rcl}
        \Sone(n,d) & \to & \Sone(n-1,d),\\
        T & \mapsto & \as_T(n) \cup \as_{\lk_T(n)}(n-1) * \{n-1\};
      \end{array}
    \right.\\
    i:
    &\left\{
      \begin{array}{rcl}
        \Sone(n-1,d) & \to & \Sone(n,d),\\
        T & \mapsto & T \cup \st_{\Hat{0}_{n,d}}(n);
      \end{array}
    \right.\\
    j:
    &\left\{
    \begin{array}{rcl}
      \Sone(n-1,d) & \to & \Sone(n,d),\\
      T & \mapsto & \lk_T(n-1) * \{n\} \cup \st_{\Hat{1}_{n,d}}(\{n-1,n\}).
    \end{array}
    \right.
  \end{align*}
\end{Lemma}

\begin{proof}
  The assertion for $f$ is contained in Theorem~\ref{thm:basics:f}.
  To prove the claims about $i$ and $j$, observe that any increasing
  flip $\Tilde{S}$ in $T \in \Sone(n-1,d)$ gives rise to an increasing
  flip $\Tilde{S}$ in $i(T)$ and an increasing flip
  $\Tilde{S} \sm \{ n-1 \} \cup \{ n \}$ in $j(T)$. This completes the proof
  of the lemma.
\end{proof}

\begin{Lemma}
  \label{thm:lemmasS1:specialsimplex}
  Let $T <_1 T' \in \Sone(n,d)$ and $S_0 := \{ n-d, \dots , n \}$.
  \begin{enumerate}\romanenumi
  \item
    If $d$ is even and $S_0$ is in $T$ then $S_0$ is also in $T'$.
  \item
    If $d$ is odd and $S_0$ is in $T'$ then $S_0$ is also in $T$.
  \end{enumerate}
\end{Lemma}

\begin{proof}
  The claim follows from the observation that for even $d$ the
  simplex $S_0$ is an upper facet of $C(n,d+1)$, whereas for odd $d$ it
  is a lower facet of $C(n,d+1)$.
\end{proof}

\begin{Lemma}
  \label{thm:lemmasS1:composition}
  For all $T \in \Sone(n,d)$ we have $i(f(T)) \le_1 T \le_1 j(f(T))$.
\end{Lemma}

\begin{proof}
  We start with a geometric description of the flip sequences that are
  going to establish the claim.
  Think of the action of $i \circ f$ as sliding vertex $n$ of
  a triangulation $T$ of $C(n,d)$ continuously to $n-1$ along the edge
  $\{ n-1, n \}$ and then adding a collection of lower facets of
  $C(n,d+1)$ to the result. If one imagines this process taking place in
  $C(n,d+1)$ then one observes that the characteristic section of $T$
  slides to the characteristic section of $i(f(T))$. Every $d$-simplex $S$
  in $T$ that contains $n$ but not $n-1$ slides exactly across the
  $(d+1)$-simplex $S \cup \{ n-1 \}$ (see Figure~\ref{fig:slide}).
  As the characteristic section $T$
  slides, these simplices $S$ are the only ones whose paths sweep out
  $(d+1)$-dimensional simplices. This yields a set
  of $(d+1)$-simplices as in the assumptions of
  Lemma~\ref{thm:lemmasS1:S1}.

  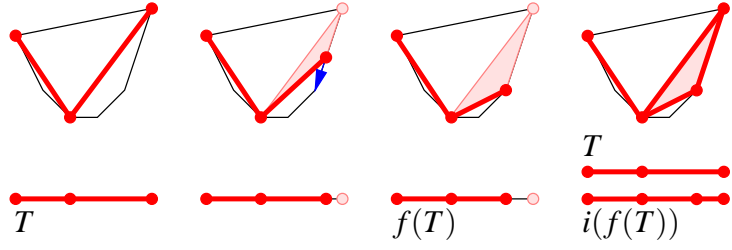
\begin{figure}[tp]
    \begin{center}
      \leavevmode
      \input{bauclic_slide.pstex_t}
      \caption{The characteristic section of $T$ slides to the
        characteristic section of $f(T)$. The simplices containing $n$ and not
        containing $n-1$ are sweeping out the increasing flips
        from $i(f(T))$ to~$T$.}
      \label{fig:slide}
    \end{center}
  \end{figure}

  On the other hand, one may regard the action of $j \circ f$ as sliding
  vertex $n-1$ of $T$ continuously to~$n$ along the edge
  $\{ n-1, n \}$ and then
  adding a bunch of upper facets fo $C(n,d+1)$ to the result.
  However, the slide---considered in $C(n,d+1)$---moves
  the characteristic section of $T$ to the characteristic
  section of $j(f(T))$. Again, the ``tracks'' of certain $d$-simplices
  provide the connecting set of $(d+1)$-simplices.
  
  In the following, we give a combinatorial
  proof of this idea.
  For $T \in \Sone(n,d)$ define the following abbreviations.
  \begin{align*}
    A(T) &:= \setof{S \in T}{n \in S, n-1 \notin S},\\
    B(T) &:= \setof{S \in T}{n \notin S, n-1 \in S}.
  \end{align*}

  We prove that $i(f(T)) \le_1 T$ for an arbitrary
  $T \in \Sone(n,d)$.
  Consider the following set of $(d+1)$-simplices in $C(n,d+1)$.
  \begin{displaymath}
    \Tilde{A}(T) := \setof{S \cup \{ n-1 \}}{S \in A(T)}.
  \end{displaymath}
  We claim that $\Tilde{A}(T)$ connects $i(f(T))$ and~$T$.  To verify
  this claim, we check properties (i)-(vi) from Lemma~\ref{thm:lemmasS1:S1}
  in Steps (i)-(vi) below.
  
  \step{(i)} All pairs of $(d+1)$-simplices in $\Tilde{A}(T)$ are admissible in
  $C(n,d+1)$ because any zig-zag-path of length $d+3$ can be transformed
  into a zig-zag-path of length $(d+2)$ by deleting $n$; deleting $n$ from
  a simplex in $\Tilde{A}(T)$, however gives a simplex in
  $f(T)$; all of these are clearly admissible in $C(n-1,d)$.

  \step{(ii)} We now show that every lower facet $S$ of a $(d+1)$-simplex
  $\Tilde{S}$ in $\Tilde{A}(T)$ is either in $i(f(T))$ or there is another
  $(d+1)$-simplex $\Tilde{S}$ in $\Tilde{A}(T)$ containing~$S$.

  To this end, let $S$ be an arbitrary lower facet of a $(d+1)$-simplex
  $\Tilde{S}$ in $\Tilde{T}$. Hence, $\Tilde{S} \sm S$ is an even gap of $S$
  in~$\Tilde{S}$.
  
  \case{1}
  If $\Tilde{S} \sm S = n$ then $S$ is contained in~$f(T)$, in particular
  it is contained in~$i(f(T))$.
  
  \case{2}
  If $\Tilde{S} \sm S = s < n-1$ then
  $F := S \sm \{n-1\}$ is a $(d-1)$-simplex in~$T$.

  If $F$ is a facet
  of $C(n,d)$ then it is an upper facet of $C(n,d)$ because $n-1$ is an
  odd gap in~$F$. Then $S = F \cup \{ n-1 \}$ was already a lower
  facet of $C(n,d+1)$
  containing $n$ and $n-1$. However, all these lower facets of $C(n,d+1)$
  are in $i(f(T))$ by construction, and thus $S \in i(f(T))$.
  
  If $F$ is not a facet of $C(n,d)$ then there is another simplex
  $S' \in T$ with $S' \neq S$ and $F \subset S'$.
  Since $n-1 \not\in S'$ we have that
  $\Tilde{S}' := S' \cup \{ n-1 \} \in \Tilde{A}(T)$
  with $\Tilde{S}' \neq \Tilde{S}$ and $S \subset \Tilde{S}'$.

  \step{(iii)} Next, we show that every upper facet $S$ of a $(d+1)$-simplex
  $\Tilde{S}$ in $\Tilde{A}(T)$ is either in $T$ or there is another
  $(d+1)$-simplex $\Tilde{S}$ in $\Tilde{A}(T)$ containing~$S$. 

  To see this, let $S$ be an arbitrary upper facet of a $(d+1)$-simplex
  $\Tilde{S}$ in $\Tilde{T}$. Hence, $\Tilde{S} \sm S$ is an odd gap of $S$ 
  in~$\Tilde{S}$.

  \case{1}
  If $\Tilde{S} \sm S = n-1$ then $S$ is contained in~$T$ by the definition
  of~$\Tilde{A}(T)$.

  \case{2}
  If $\Tilde{S} \sm S = s < n-1$ then
  $F := S \sm \{n-1\}$ is a $(d-1)$-simplex in~$T$.

  We show now that $F$ is not a facet of $C(n,d)$:
  Because $s$ is an odd gap of $S$ in $\Tilde{S}$ and $n-1 > s$ is an
  additional gap of $F$ larger than $s$ we conclude that
  $s$ is an even gap of $F$ in $\Tilde{S}$. However, $n-1$ is clearly
  an odd gap of $F$ in $\Tilde{S}$ because $n \in F$. Thus, $F$ contains
  an even and an odd gap, and is therefore not a facet of~$\Tilde{S}$.
  Consequently, it cannot be a facet of $C(n,d)$.
 
  Hence, there is another simplex
  $S' \in T$ with $S' \neq S$ and $F \subset S'$.
  Since $n-1 \not\in S'$ we have that
  $\Tilde{S}' := S' \cup \{ n-1 \} \in \Tilde{A}(T)$
  with $\Tilde{S}' \neq \Tilde{S}$ and $S \subset \Tilde{S}'$.

  \step{(iv)}
  We now prove that every $d$-simplex in $i(f(T)) \sm T$ is a
  lower facet of some $(d+1)$-simplex $\Tilde{S}$ in~$\Tilde{A}(T)$.

  Let $S$ be a $d$-simplex in $i(f(T))$ but not in $T$.
  There are two types of $d$-simplices in $i(f(T)) \sm T$: simplices
  of the form $S = S' \sm \{n\} \cup \{n-1\}$ with $S' \in A(T)$ (case~1), and
  lower facets of $C(n,d+1)$ containing $n$ and $n-1$ (case~2) .

  \case{1}
  If $S$ is of the form $S = S' \sm \{n\} \cup \{n-1\}$ with $S' \in A(T)$
  then $\Tilde{S} := S \cup \{ n \}$ is in $\Tilde{A}(T)$, and
  $n$ is clearly an even gap of $S$ in $\Tilde{S}$. Thus,
  $S$ is a lower facet of the simplex $\Tilde{S} \in \Tilde{A}(T)$.

  \case{2}
  If $S$ is a lower facet of $C(n,d+1)$ containing $n$ and $n-1$ then
  all gaps of $S$ are even. Hence, all gaps of
  $F := S \sm \{n-1\}$ are odd. Thus, $F$ is an upper facet of
  $C(n,d)$. This leads to the existence of a $d$-simplex $S'$ in $T$
  containing~$F$. Since $n$ is in $F$ we know that $n$ is also in $S'$.
  If $n-1 \in S'$ then
  $S = S' \in T$; contradiction to $S \in i(f(T)) \sm T$.
  Therefore, $S'$ is in $A(T)$ and, consequently,
  $\Tilde{S} := S' \cup \{n-1\}$ is a $(d+1)$-simplex in $\Tilde{A}(T)$.
  Moreover, $S = F \cup \{ n-1 \}$ is a facet of~$\Tilde{S}$ because
  $\Tilde{S}$ contains $n-1$. Additionally, $S$ is ---
  by the assumption of this case ---
  a lower facet of $C(n,d+1)$, so it must be a lower facet of~$\Tilde{S}$.

  \step{(v)}
  We now prove that every $d$-simplex in $T \sm i(f(T))$ is an
  upper facet of some $(d+1)$-simplex $\Tilde{S}$ in~$\Tilde{A}(T)$.

  Let $S$ be a $d$-simplex in $T$ but not in $i(f(T))$.
  Then $S$ is, in particular, not contained in $f(T)$. There are two
  types of $d$-simplices in $T \sm f(T)$: simplices from $A(T)$ (case~1),
  and simplices containing both $n-1$ and~$n$ (case~2).

  \case{1}
  Assume $S$ is in $A(T)$. Then $\Tilde{S} := S \cup \{ n-1 \}$ is
  in $\Tilde{A}(T)$.
  Since $n-1$ is an odd gap of $S$ in $\Tilde{S}$ we conclude that
  $S$ is an upper facet of the $(d+1)$-simplex $\Tilde{S}$ in~$\Tilde{A}(T)$.

  \case{2}
  If both $n-1$ and $n$ are in~$S$ then $S$ cannot be a lower facet of
  $C(n,d+1)$, because all lower facets of $C(n,d+1)$ containing
  both $n$ and $n-1$ are in $i(f(T))$ by construction.
  Assume, for the sake of contradiction, that $F := S \sm \{n-1\}$
  is a facet of $C(n,d)$. Then either all gaps of
  $F$ are even or all gaps of $F$ are odd. Since $n \in F$ we know that
  $n-1$ is an odd gap of $F$, thus all gaps of $F$ must be odd. However,
  then all gaps of $S = F \cup \{n-1\}$
  are even; contradiction to the fact that
  $S$ is not a lower facet of $C(n,d+1)$. We conclude that
  $F$ is not a facet of $C(n,d)$. Thus, there is another simplex
  $S' \neq S$ in $T$ containing~$F$. Moreover, because $n-1 \not\in S'$
  but $n \in S'$, we have $S' \in A(T)$, and, consequently,
  $\Tilde{S} := S' \cup \{ n-1 \}$ is in $\Tilde{A}(T)$.
  Because $n-1$ is an odd gap of $S'$ in $\Tilde{S}$ we know that
  $S'$ is an upper facet of $\Tilde{S}$.
  Moreover, since $S'$ and $S$ are both in $T$ they are admissible in
  $C(n,d)$. That means that $S$ is also an upper facet of $\Tilde{S}$.
  (A lower and an upper facet of a $(d+1)$-simplex in $C(n,d+1)$
  are never admissible in $C(n,d)$.)
  
  \step{(vi)}
  Finally we prove that every simplex in $T \sm i(f(T)) \cup i(f(T)) \sm T$
  is a facet of at most one $(d+1)$-simplex in $\Tilde{A}(T)$.

  \case{1} $S$ is a $d$-simplex in $T \sm i(f(T))$. If $S$ is in $A(T)$ then
  $\Tilde{S} = S \cup \{n-1\}$ is the only $(d+1)$-simplex in $\Tilde{A}(T)$
  containing~$S$ because membership in $\Tilde{A}(T)$ requires the containment
  of $n-1$. If both $n$ and $n-1$
  are in $S$ then we proceed as follows. Assume, for the sake of
  contradiction, that there are
  two distinct $(d+1)$-simplices $\Tilde{S}$ and $\Tilde{S}'$ in
  $\Tilde{A}(T)$ containing~$S$. Then $S$ is a lower facet of one of them,
  say $\Tilde{S}$ and an upper facet of the other one, say $\Tilde{S}'$.
  In other words, $s := \Tilde{S} \sm \Tilde{S}' < n-1$ is an odd gap of
  $\Tilde{S}'$ and $s' := \Tilde{S}' \sm \Tilde{S} < n-1$ is an even gap
  of $\Tilde{S}$ in $\Tilde{S} \cup \Tilde{S}'$. 
  By construction of
  $\Tilde{A}(T)$, we know that $\Tilde{S} = R \cup \{n-1\}$
  and $\Tilde{S}' = R' \cup \{n-1\}$ for some $R, R' \in A(T)$.
  In particular, $R$ and $R'$ are in $T$, thus admissible in $C(n,d)$.
  However, $R = \Tilde{S} \sm \{n-1\}$ and $R' = \Tilde{S}' \sm \{n-1\}$.
  Therefore, $s$ is an even gap of $R'$ and $s'$ is an odd gap of $R$ in
  $R \cup R'$. But that means, $R'$ is a lower and $R$ is an upper facet
  of the $(d+1)$-simplex $R \cup R'$; contradiction to the fact that
  $R$ and $R'$ are admissible in $C(n,d)$.

  \case{2} $S$ is a $d$-simplex in $i(f(T)) \sm T$. If $S$ is of
  the form $S' \sm \{n\} \cup \{n-1\}$ for some $S' \in A(T)$ then
  $S \cup \{n\}$ is the only $(d+1)$-simplex in $\Tilde{A}(T)$ containing~$S$
  because membership in $\Tilde{A}(T)$ requires the containment of~$n$.
  On the other hand, if $S$ is a lower facet of $C(n,d+1)$
  then there cannot be two distinct
  $(d+1)$-simplices which both contain $S$ and are admissible in $C(n,d+1)$.
  
  Steps 1 to~6 prove that
  the assumptions of Lemma~\ref{thm:lemmasS1:S1} are satisfied, thus
  $\Tilde{A}(T)$ connects $i(f(T))$ and~$T$, proving $i(f(T)) \le_1 T$.

  Analogously, the set
  \begin{displaymath}
    \Tilde{B}(T) := \setof{\Tilde{S} \cup \{ n \}}{S \in B(T)}.
  \end{displaymath}
  connects $T$ and $j(f(T))$, proving $T \le_1 j(f(T))$.  We omit
  the details verifying this, which are similarly tedious.
\end{proof}

\begin{Lemma}
  \label{thm:lemmasS1:fibers01}
  Let $T$ be in $\Sone(n,d)$ and $S_0 := \{ n-d, \dots , n \}$.
  \begin{enumerate}\romanenumi
  \item \label{itm:fibers01:even0}
    Let $d$ be even, $f(T) = \Hat{0}_{n-1,d}$, and $S_0 \notin T$.
    Then $T = \Hat{0}_{n,d}$.
  \item \label{itm:fibers01:even1}
    Let $d$ be even, $f(T) = \Hat{1}_{n-1,d}$, and $S_0 \in T$.
    Then $T = \Hat{1}_{n,d}$.
  \item \label{itm:fibers01:odd0}
    Let $d$ be odd, $f(T) = \Hat{0}_{n-1,d}$, and $S_0 \in T$.
    Then $T = \Hat{0}_{n,d}$.
  \item \label{itm:fibers01:odd1}
    Let $d$ be odd, $f(T) = \Hat{1}_{n-1,d}$, and $S_0 \notin T$.
    Then $T = \Hat{1}_{n,d}$.
  \end{enumerate}
\end{Lemma}

\begin{proof}
  For the proof of \ref{itm:fibers01:even0}, let $T \in \Sone(n,d)$ for
  even~$d$ with $f(T) = \Hat{0}_{n-1,d}$. Assume that $T \neq
  \Hat{0}_{n,d}$. Recall that any such element $T$ in $\Sone(n,d)$ can be
  connected to $\Hat{0}_{n,d}$ by a sequence of decreasing
  flips (see Theorem~\ref{thm:basics:bounded}).
  The map $f$ is order-preserving (see Theorem~\ref{thm:basics:f}); thus every
  element in such a sequence is mapped by $f$ to~$\Hat{0}_{n-1,d}$.
  Because of Lemma~\ref{thm:lemmasS1:specialsimplex}, we may therefore
  assume that $T$ differs from $\Hat{0}_{n,d}$ by exactly one increasing
  flip corresponding to a $(d+1)$-simplex~$\Tilde{S}$. The simplex
  $\Tilde{S}$ must contain both $n-1$ and~$n$ because otherwise it would
  give rise to a (non-trivial)
  flip from $\Hat{0}_{n-1,d}$ to $f(T)$ in contradiction to
  $f(T) = \Hat{0}_{n-1,d}$.
  The only flip in $\Hat{0}_{n,d}$ containing $n-1$ and~$n$ corresponds,
  however, to the $(d+1)$-simplex $\Tilde{S} = \{ 1, n-d, \dots , n \}$.
  The fact that $S_0$ is an upper facet of $\Tilde{S}$, thus is contained
  in the triangulation resulting from this flip, gives a
  contradiction.
  Thus, claim~\ref{itm:fibers01:even0} is proved.

  The proofs of the remaining statements are analogous with
  \begin{displaymath}
    \Tilde{S} =
    \begin{cases}
      \{ n-d-1,n-d, \dots , n \} &
      \text{decreasing flip in $\Hat{1}_{n,d}$ for~\ref{itm:fibers01:even1}},\\
      \{ n-d-1,n-d, \dots , n \} &
      \text{increasing flip in $\Hat{0}_{n,d}$ for~\ref{itm:fibers01:odd0}},\\
      \{ 1,n-d, \dots , n \} &
      \text{decreasing flip in $\Hat{1}_{n,d}$ for~\ref{itm:fibers01:odd1}}.
    \end{cases}
  \end{displaymath}
\end{proof}

%%%%%%%%%%%%%%%%%%%%%%%%%%%%%%%%%%%%%%%%%%%%%%%%%%%%%%%%%%%%%%%%%%%%%%%%%%%%%%%
%% Lemmas on $\Stwo(n,d)$:                                                   %%
%% 07/03/97 pe/jr/vr                                                         %%
%%%%%%%%%%%%%%%%%%%%%%%%%%%%%%%%%%%%%%%%%%%%%%%%%%%%%%%%%%%%%%%%%%%%%%%%%%%%%%%

\section{Lemmas for $\Stwo(n,d)$}
\label{sec:lemmasS2}

This section is devoted to proving an analogous set of lemmas to
the ones in the previous section, in order to guarantee the assumptions of the
Suspension Lemma for $\Stwo(n,d)$.
Again, in the following $n > d+2$.

\begin{Lemma}
  \label{thm:lemmasS2:maps}
  The following maps are order-preserving.
  \begin{align*}
    f:
    &\left\{
      \begin{array}{rcl}
        \Stwo(n,d) & \to & \Stwo(n-1,d),\\
        T & \mapsto & T \sm n := \as_T(n) \cup \as_{\lk_T(n)}(n-1) * \{n-1\};
      \end{array}
    \right.\\
    i:
    &\left\{
      \begin{array}{rcl}
        \Stwo(n-1,d) & \to & \Stwo(n,d),\\
        T & \mapsto & T \cup \st_{\Hat{0}_{n,d}}(n);
      \end{array}
    \right.\\
    j:
    &\left\{
    \begin{array}{rcl}
      \Stwo(n-1,d) & \to & \Stwo(n,d),\\
      T & \mapsto & \lk_T(n-1) * \{n\} \cup \st_{\Hat{1}_{n,d}}(\{n-1,n\}).
    \end{array}
    \right.
  \end{align*}
\end{Lemma}

\begin{proof}
  That $i$ and $j$ are order-preserving is easily seen by considering
  the following facts: both maps embed a triangulation of $C(n-1,d)$ into
  $C(n,d)$; $i$ copies the original triangulation, $j$ renames
  $n-1$ to~$n$. This does not change any height relations of piecewise
  linear sections to each other.
  Then both maps add a set of simplices which does not depend
  upon~$T$. These are consequently at the same height for all
  triangulations. Thus, all height relations are maintained.

 We now prove the assertion concerning~$f$.   We use the fact that
 the map $f: \Sonetwo(n,d) \rightarrow \Sonetwo(n-1,d)$ has the
 following geometric interpretation: given a triangulation $T$ of
 $C(n,d)$, imagine a homotopy that ``slides" the vertex $n$ down
 the moment curve toward the vertex $n-1$, so that at $t=0$ one
 has the triangulation $T(0)=T$ of the original cyclic polytope $C(n,d)$,
 and at $t=1$ some of the simplices of $T(1)$ (namely those containing
 both $n-1$ and $n$) have become degenerate (volume zero).  If one
 eliminates these degenerate simplices from $T(1)$ and relabels the
 vertex $n$ by $n-1$ in the remaining simplices, one obtains the triangulation
 $f(T)$ of $C(n-1,d)$.

 To prove that $f$ is order-preserving, assume $T \leq_2 T'$, and
 we will show that $f(T) \leq_2 f(T').$  Fix a point $x \in C(n-1,d)$,
 and for $0 \leq t \leq 1$, let $T(t)(x)_{d+1}, T'(t)(x)_{d+1}$ be the
 $(d+1)^{st}$-coordinates  of the image of $x$ under the
 parametrized characteristic sections 
 $T(t), T'(t): C(n,d) \rightarrow C(n,d+1)$.
 Since $T \leq_2 T'$, we have
 
\begin{displaymath}
T'(t)(x)_{d+1} - T(t)(x)_{d+1} \geq 0 \,\, \text{ for }\,\, 0 \leq t < 1.
\end{displaymath}

\noindent
 However $T'(t)(x)_{d+1} - T(t)(x)_{d+1}$ is clearly a continuous
 function of $t$, so the same inequality holds for $t=1$.  Hence
  
\begin{displaymath}
f(T)(x)_{d+1}= T(1)(x)_{d+1} \leq T'(1)(x)_{d+1} = f(T')(x)_{d+1}
\end{displaymath}

\noindent
 which shows that $f(T) \leq_2 f(T')$

\end{proof}

\begin{Lemma}
  \label{thm:lemmasS2:specialsimplex}
  Let $T < T' \in \Stwo(n,d)$ and $S_0 := (n-d, \dots , n)$.
  \begin{enumerate}\romanenumi
  \item
    If $d$ is even and $S_0$ is in $T$ then $S_0$ is also in $T'$.
  \item
    If $d$ is odd and $S_0$ is in $T'$ then $S_0$ is also in $T$.
  \end{enumerate}
\end{Lemma}

\begin{proof}
  The assertion follows from exactly the same argument as given
  in the proof of Lemma~\ref{thm:lemmasS1:specialsimplex}.
\end{proof}

\begin{Lemma}
  \label{thm:lemmasS2:composition}
  For all $T \in \Stwo(n,d)$ we have $i(f(T)) \le_2 T \le_2 j(f(T))$.
\end{Lemma}

\begin{proof}
  This follows from Lemma~\ref{thm:lemmasS1:composition} and the fact
  that $T \le_1 T'$ always implies $T \le_2 T'$
  (see~\cite{EdelmanReiner1996}).
\end{proof}

\begin{Lemma}
  \label{thm:lemmasS2:fibers01}
  Let $T$ be in $\Stwo(n,d)$ and $S_0 := \{ n-d, \dots , n \}$.
  \begin{enumerate}\romanenumi
  \item
    Let $d$ be even, $f(T) = \Hat{0}_{n-1,d}$, and $S_0 \notin T$.
    Then $T = \Hat{0}_{n,d}$.
  \item
    Let $d$ be even, $f(T) = \Hat{1}_{n-1,d}$, and $S_0 \in T$.
    Then $T = \Hat{1}_{n,d}$.
  \item
    Let $d$ be odd, $f(T) = \Hat{0}_{n-1,d}$, and $S_0 \in T$.
    Then $T = \Hat{0}_{n,d}$.
  \item
    Let $d$ be odd, $f(T) = \Hat{1}_{n-1,d}$, and $S_0 \notin T$.
    Then $T = \Hat{1}_{n,d}$.
  \end{enumerate}
\end{Lemma}

\begin{proof}
  This statement is independent of the partial order
  $\Sone(n,d)$ or $\Stwo(n,d)$ under consideration. Thus the
  proof of Lemma~\ref{thm:lemmasS1:fibers01} is valid here as well.
\end{proof}

%%%%%%%%%%%%%%%%%%%%%%%%%%%%%%%%%%%%%%%%%%%%%%%%%%%%%%%%%%%%%%%%%%%%%%%%%%%%%%%
%% The generalized Baues problem for $C(n,d)$ with $d \leq 3$:               %%
%% 07/03/97 pe/jr/vr                                                         %%
%%%%%%%%%%%%%%%%%%%%%%%%%%%%%%%%%%%%%%%%%%%%%%%%%%%%%%%%%%%%%%%%%%%%%%%%%%%%%%%

\section{The Generalized Baues Problem for $C(n,d)$ with $d \leq 3$}
\label{sec:baues}

The goal of this section is to prove a new special case of the
generalized Baues problem,
but we must first recall the definition of the Baues poset $\oo$.
A \emph{polytopal decomposition} $\delta$ of $C(n,d)$ is a collection
$\{V_\alpha\}$ of vertex subsets $V_\alpha \subseteq [n]$ satisfying
\begin{itemize}
\item For all $\alpha$, $|V_\alpha| \geq d+1$.
\item Any two cyclic subpolytopes $C(V_\alpha,d), C(V_\beta,d)$ intersect
in a common face (possibly empty).
\item The union of the cyclic subpolytopes $C(V_\alpha,d)$ covers
  $C(n,d)$, i.~e.,
  \begin{displaymath}
    \bigcup_\alpha C(V_\alpha,d) = C(n,d)
  \end{displaymath} 
\end{itemize}
Say that a polytopal decomposition is \emph{proper} if it is
not the trivial decomposition $\{[n]\}$. 

The \emph{Baues poset} $\oo$ is the set of all proper
polytopal decompositions ordered by refinement, i.~e.,
$\delta = \{V_\alpha\} \leq \delta' = \{V_{\alpha'}\}$
if for every $V_\alpha \in \delta$ there exists a
$V_{\alpha'} \in \delta'$ with 
$V_\alpha \subseteq V_{\alpha'}$.  One can check that this
agrees with the poset considered in the
Generalized Baues Problem \cite{BilleraKapranovSturmfels1994}
for the case of subdivisions of a cyclic polytope.
Theorem~\ref{thm:introduction:baues3main} now reads as follows.

\begin{Theorem}
\label{thm:main-Baues-thm}
  For $d \le 3$ the poset $\oo$ is homotopy equivalent
  to a sphere of dimension $n-d-2$.
\end{Theorem}

As was said in the introduction, our method will be to show that
$\oo$ is homotopy equivalent to the suspension $\susp(\proper{\Stwo(n,d)})$.
We begin by defining a map $\phi$ from $\oo$ to intervals in $\Stwo(n,d)$.
An element $\delta$ of $\oo$ is a polytopal
subdivision of $C(n,d)$, so let $\phi(\delta)$ be the set of all triangulations
of $C(n,d)$ which refine it.

\begin{Lemma}
\label{lemma:Baues:define-map}
For any $\delta$ in $\oo$,
\begin{itemize}
\item the set $\phi(\delta)$ is a non-empty interval in $\Stwo(n,d)$.
\item $\phi(\delta)$ is \emph{not} the improper interval consisting of
all $\Stwo(n,d)$.
\item $\delta \leq \delta'$ in $\oo$ implies 
$\phi(\delta) \subseteq \phi(\delta')$.
\item $\phi$ is injective, i.~e., $\phi(\delta)=\phi(\delta')$ implies 
$\delta=\delta'$.
\end{itemize}
\end{Lemma}
\begin{proof}
Since $\delta$ is a polytopal subdivision of $C(n,d)$, and
subsets $V$ of the vertices of $C(n,d)$ span cyclic subpolytopes
$C(V,d)$, we know that $\delta$ gives a decomposition
\begin{displaymath}
  C(n,d) = \bigcup_{\alpha} C(V_\alpha, d)
\end{displaymath}
for some vertex sets $V_\alpha$ in which the $C(V_\alpha,d)$ all
have disjoint interiors.
If we let $\hat{0}_\alpha,\hat{1}_\alpha$ denote the bottom
and top triangulations of $C(V_\alpha, d)$, then one can
form two triangulations $T$ and $T'$ respectively,
by refining $\delta$ according to $\hat{0}_\alpha$ and $\hat{1}_\alpha$
respectively on each subpolytope $C(V_\alpha,d)$.
It is then clear from the definition of $\Stwo(n,d)$ that $\phi(\delta) = [T,T']$.
This proves the first assertion of the lemma.
 
To prove the second assertion, note that since $\delta$ is
a non-trivial polytopal subdivision of $C(n,d)$, it must
use at least one $(d-1)$-simplex $\sigma$ spanned by the vertices of
$C(n,d)$ which lies interior to $C(n,d)$, and therefore this
simplex $\sigma$ would lie in every triangulation in $\phi(\delta)$.
If $\phi(\delta)$ were all of $\Stwo(n,d)$, then in particular this would
imply that the bottom and top triangulations $\hat{0}, \hat{1}$
have this simplex $\sigma$ in common.  But one can easily
check from the explicit description of the triangulations 
$\hat{0}, \hat{1}$ given in \cite{EdelmanReiner1996} or \cite{Rambau1996}
that they have no interior $(d-1)$-simplices in common.

To see the third assertion, note $\delta \leq \delta'$ means that
$\delta$ refines $\delta'$ as a polytopal subdivision,
so any triangulation $T$  which refines $\delta$ 
will also refine $\delta'$, and hence $\phi(\delta) \subseteq \phi(\delta')$.

To see the fourth assertion, it suffices to show that 
$\delta$ is completely determined by $\phi(\delta)$, in the
sense that the set of $(d-1)$-simplices of $\delta$ is the
intersection of all the sets of $(d-1)$-simplices of its triangulation
refinements.  Certainly the $(d-1)$-simplices
of $\delta$ are contained in this intersection.
This intersection cannot be larger because for each $\alpha$,
(using the notation of the first paragraph),
the two triangulations $\hat{0}_\alpha, \hat{1}_\alpha$
share no common $(d-1)$-simplices interior to $C(V_\alpha,d)$. 
\end{proof}

We next recall and introduce some notions about lattices.
Given a lattice $L$ with bottom and top elements $\hat{0}, \hat{1}$,
an element of $L$ which covers $\hat{0}$
(resp. is covered by $\hat{1}$) is called an \emph{atom (coatom)}, resp.
The lattice $L$ is \emph{atomic} (resp.\ \emph{coatomic})
if the join of all the atoms
is $\hat{1}$ (resp. the meet of all the coatoms is $\hat{0}$).
Any interval $[x,y]$ in a lattice is a lattice itself, and will be
called atomic or coatomic if it satisfies the previous conditions.
An interval $[x,y]$ will be called \emph{proper} if it is not the
whole lattice $L=[\hat{0},\hat{1}]$.  Recall that
the \emph{proper part} of $L$ is the subposet 
$\proper{L}: = L \sm \{\hat{0},\hat{1}\}$.

We now define three \emph{interval posets} 
as certain collections of intervals in $L$ ordered by inclusion of
intervals:
\begin{itemize}
\item $\Int(L)$ --- all non-empty intervals in $L$,
\item $\intbar(L)$ --- all non-empty, proper intervals in $L$,
\item $\intbar_{\atomic}(L)$ --- all non-empty, proper, atomic intervals in $L$.
\end{itemize}
Similarly one can define $\intbar_{\coatomic}(L)$.

In \cite{Walker1981} it was shown that $\Int(L)$ is canonically
homeomorphic to $L$, and that $\intbar(L)$ is canonically homeomorphic
to $\susp(\proper{L})$, i.~e., the suspension of the proper part of $L$.
One can view Lemma~\ref{lemma:Baues:atomic-proper-intervals} below
as asserting an analogous statement,
up to homotopy, for $\intbar_{\atomic}(L)$.

We recall (Theorem~\ref{latticeness})  that for $d \leq 3$
$\Stwo(n,d)$ is a lattice, and note that Lemma~\ref{lemma:Baues:define-map}
shows that $\phi$ defines an injective, order-preserving map
$\oo \rightarrow \intbar(\Stwo(n,d))$.
  
\begin{Lemma}
For $d \leq 3$, the image of $\phi:\oo \rightarrow \intbar\Stwo(n,d)$ is
exactly $$\intbar_{\coatomic}(\Stwo(n,d)).$$
\end{Lemma}
\begin{proof}
To see that $\phi(\delta)$ is always a coatomic interval in $\Stwo(n,d)$,
we use the notation from the proof of Lemma~\ref{lemma:Baues:define-map},
and note the following isomorphism of posets:
\begin{displaymath}
  \phi(\delta) = [T,T'] \cong \prod_{\alpha} [\hat{0}_\alpha,\hat{1}_\alpha].
\end{displaymath}
Since each interval $[\hat{0}_\alpha,\hat{1}_\alpha]$ is
isomorphic to $\Stwo(n',d)$ for some $n'<n$, the coatomicity of $\phi(\delta)$
would follow if we knew that $\Stwo(n,d)$ is a coatomic
lattice for $d \leq 3$.  But if $\Stwo(n,d)$ were \emph{not} coatomic
then its proper part $\proper{\Stwo(n,d)}$ would be contractible
(see, e.~g.,\cite[Theorem~10.14]{Bjoerner1995}),
contradicting Theorem~\ref{Tamari-sphericity} above.

It remains then to show that every coatomic interval in $\Stwo(n,d)$
is of the form $\phi(\delta)$ for some $\delta$ in $\oo$.
For $d=1$, this is trivial since the cyclic polytope $C(n,1)$
is simply a line segment with $n-2$ interior subdivision points.
Triangulations of $C(n,1)$
are specified by their subset of interior vertices
and $\Stwo(n,d)$ is a Boolean algebra $\mathcal B_{n-2}$,
so that every interval is coatomic, and it is easy to see that
every interval is $\phi(\delta)$ for some $\delta$ in $\oo$.
 
For $d=2, 3$ the fact that every coatomic interval in $\Stwo(n,d)$
is of the form $\phi(\delta)$ requires some argument.  Assume
we have such a coatomic interval $[T,T']$, and we will show
how to construct its preimage $\delta$.  Form a graph $G$
whose vertices are the $d$-simplices $\sigma$ in the triangulation
$T'$, and whose edges correspond to a pair of $d$-simplices
$\sigma, \sigma'$ which share a $(d-1)$-simplex $\tau$ that is \emph{not}
a simplex in $T$.  Let $\{G_\alpha\}$ be the various connected components
of $G$, and define $V_\alpha$ to be the set of all vertices of $C(n,d)$
which lie in a simplex of $G_\alpha$.  We wish to prove two
claims about these graphs:
\begin{itemize}
\item If $\sigma, \sigma'$ are simplices of $T'$ which correspond
to an edge of $G$, then their union is a cyclic subpolytope $C(d+2,d)$
which supports a bistellar operation corresponding to a covering
relation between $T'$ and some coatom  of the interval $[T,T']$.
\item For each $\alpha$, the connected component $G_\alpha$ is a path, and
the set of $d$-simplices $\sigma$ corresponding to $G_\alpha$
gives exactly the maximal simplices of
the top triangulation $\hat{1}_\alpha$ of the
cyclic subpolytope $C(V_\alpha,d)$.
\end{itemize}
Assuming these two claims for the moment, we show how to finish the proof.
The second claim implies that the decomposition 
$C(n,d) = \bigcup_\alpha C(V_\alpha,d)$ defines a polytopal subdivision
$\delta$.  Furthermore, as in the first paragraph of this proof,
we know that  $\phi(\delta)$ is equal to some coatomic interval 
$[T_\delta, T'_\delta]$, where $T,T'$ refine $\delta$ and the
restriction to $C(V_\alpha,d)$ of $T, T'$
looks like $\hat{0}_\alpha, \hat{1}_\alpha$ respectively.
By the second claim, this means that $T' = T'_\delta$.
By both claims together, every coatom of the interval 
$[T_\delta,T'_\delta]$ is also a coatom of $[T,T']$
(i.~e., all of the former coatoms lie above $T$), and hence by coatomicity
of $[T,T']$ we must have $T=T_\delta$.
Therefore $[T,T'] = \phi(\delta)$ as desired.

To show the first claim, assume $\sigma, \sigma'$ are
simplices of $T'$ which correspond to an edge of $G$, so
there intersection is a $(d-1)$-simplex $\tau$ which is not in $T$.
Assume for the sake of contradiction that
the union $\sigma \cup \sigma'$ does \emph{not} support
a bistellar operation as asserted in the claim.  Then every
coatom $T''$ of $[T,T']$ will have $\tau$ in its \emph{submersion set}
$\submersion(T'')$ (see Proposition~\ref{submersion}).
Since the meet operation in $\Stwo(n,d)$ corresponds to intersection of
submersion sets, coatomicity of $[T,T']$ implies that
$\submersion(T)$ would also contain $\tau$.  But then the fact that
$\tau$ is not a $(d-1)$-simplex of $T$ would imply that
\begin{itemize}
\item if $d=2$ then $\tau = \{i,j\}$ must be ``foiled" by some other
$\tau' = \{k,l\}$ in $\submersion(T)$ which satisfies $i < k < j < l$
(see Proposition~\ref{low-dimensional-submersion}).
\item if $d=3$ then $\tau = \{i,j,k\}$ must be ``foiled" by
one of its edges, say $\{i,j\}$, \emph{intertwining}
another triple $\tau' = \{x,y,z\}$ in $\submersion(T)$
in the sense that $x < i < y < j < z$
(see Proposition~\ref{low-dimensional-submersion})
\end{itemize}
However in both of these cases, $\tau'$ would also lie in
$\submersion(T')$ since $T < T'$ in $\Stwo(n,d)$, and hence would ``foil"
$\tau$ from being a $(d-1)$-simplex of $T'$. Contradiction.

To show the second claim, note that the first claim implies
very stringent requirements on what $\sigma, \sigma'$
can look like whenever they correspond to an edge in $G$:
\begin{itemize}
\item if $d=2$, $\sigma = \{i,j,l\}, \sigma'=\{j,k,l\}$
for some $i < j < k < l$, and
\item if $d=3$, $\sigma = \{i,j,k,m\}, \sigma'=\{i,k,l,m\}$
for some $i < j < k < l < m$.
\end{itemize}
It is easy to check that these requirements, combined with the fact
that a $(d-1)$-simplex $\tau$ can lie in at most two $d$-simplices of $T'$,
implies that the degree of any vertex in a connected component $G_\alpha$
can be at most $2$.  In fact, $G_\alpha$ is constrained
to look like the following path of $d$-simplices:
\begin{itemize}
\item for $d=2$, 
\begin{displaymath}
  \{v_1 v_2 v_r\},\{ v_2 v_3 v_r\}, \{v_3 v_4 v_r\},
  \ldots, \{v_{r-2} v_{r-1} v_r\} 
\end{displaymath}
\item for $d=3$, 
  \begin{displaymath}
    \{v_1 v_2 v_3 v_r\}, \{v_1 v_3 v_4 v_r\},\{ v_1 v_4 v_5 v_r\}, \ldots,
    \{v_1 v_{r-2} v_{r-1} v_r\}
  \end{displaymath}
\end{itemize}
where $v_1 < \cdots < v_r$ are the vertices $V_\alpha$ of $G_\alpha$
written increasing order.
In both cases this description matches exactly the top triangulation
$\hat{1}_\alpha$ of $C(V_\alpha,d)$.
\end{proof}

Once the image of $\phi$ has been established,
Theorem \ref{thm:main-Baues-thm} follows by combining 
\begin{itemize}
\item Lemma~\ref{lemma:Baues:atomic-proper-intervals} below,
\item the above-mentioned fact that the proper interval poset
$\intbar(L)$  is homeomorphic to $\susp(\proper{L})$, and
\item Theorem~\ref{Tamari-sphericity} or
  Theorem~\ref{thm:introduction:hstmain}. 
\end{itemize}

\begin{Lemma}
\label{lemma:Baues:atomic-proper-intervals}
Let $L$ be any finite lattice.
Then $\intbar_{\atomic}(L)$ (or $\intbar_{\coatomic}(L)$) is homotopy
equivalent to $\intbar(L)$.
\end{Lemma}

Lemma~\ref{lemma:Baues:atomic-proper-intervals}
follows from a more general lemma, which we think is of
independent interest.  We are indebted to P. Webb for the
statement and proof of this lemma.

\begin{Lemma}
\label{Webb-lemma}
Let $P$ be a poset with $\hat{0},\hat{1}$.  If
$\{[x_i,y_i]\}_{i=1}^r$ is any finite collection of
intervals with the open intervals $(x_i,y_i)$ contractible
for all $i$, then the inclusion
\begin{displaymath}
  \intbar P  \sm \{[x_i,y_i]\}_{i=1}^r \hookrightarrow \intbar P
\end{displaymath}
induces a homotopy equivalence.
\end{Lemma}

Lemma~\ref{lemma:Baues:atomic-proper-intervals} then follows from 
Lemma~\ref{Webb-lemma} by letting $P=L$ and letting
$\{[x_i,y_i]\}_{i=1}^r$ be the non-coatomic intervals of $L$.
These non-coatomic intervals satisfy the hypothesis of the
lemma by \cite[Theorem 10.14]{Bjoerner1995}.

Lemma~\ref{Webb-lemma} follows immediately from the
following two sublemmas:

\begin{Sublemma}\cite{Bouc1984}
\label{Bouc-lemma}
In a poset $Q$, if $\{q_i\}_{i=1}^r$ is a finite subset
with $Q_{<q_i}$ contractible for all $i$, then the inclusion
\begin{displaymath}
  Q \sm \{q_i\}_{i=1}^r \hookrightarrow Q
\end{displaymath}
induces a homotopy equivalence.
\end{Sublemma}
\begin{proof}
Re-index the elements $\{q_i\}_{i=1}^r$
in such a way that $q_i > q_j$ in $Q$ implies $i<j$.
Then 
\begin{displaymath}
  (Q \sm \{q_1,\ldots,q_{i-1}\})_{< q_i} = Q_{<q_i}
\end{displaymath}
is contractible for all $i$, so an application of
Quillen's Fiber Lemma \cite[Theorem 10.5]{Bjoerner1995}
proves the homotopy equivalence by induction on $i$.
\end{proof}

We can apply Sublemma~\ref{Bouc-lemma} with $Q=\intbar P$ to 
prove Lemma~\ref{Webb-lemma} once we have established

\begin{Sublemma}
In a poset $P$ with $\hat{0},\hat{1}$, 
if an open interval $(x,y)$ is contractible, then
$\left( \intbar P \right)_{<[x,y]}$ is contractible.
\end{Sublemma}
\begin{proof}
Note that
\begin{displaymath}
  \left( \intbar P \right)_{<[x,y]}  = \intbar [x,y].
\end{displaymath}
But $\intbar [x,y]$ is homeomorphic to
the suspension $\susp (x,y)$ by \cite{Walker1981},
and hence contractible since $(x,y)$ was assumed contractible.
\end{proof}

%%%%%%%%%%%%%%%%%%%%%%%%%%%%%%%%%%%%%%%%%%%%%%%%%%%%%%%%%%%%%%%%%%%%%%%%%%%%%%%
%% Open problems:                                                            %%
%% 24/03/97 vr                                                               %%
%%%%%%%%%%%%%%%%%%%%%%%%%%%%%%%%%%%%%%%%%%%%%%%%%%%%%%%%%%%%%%%%%%%%%%%%%%%%%%%

\section{Open Problems}
\label{sec:problems}
The following are some remarks and remaining open problems about
triangulations of cyclic polytopes which we consider interesting.
\begin{enumerate}
\item The proof of Theorem~\ref{thm:main-Baues-thm}
relied heavily on the fact established in \cite{EdelmanReiner1996}
that $\Stwo(n,d)$ is a lattice for $d \leq 3$.
Unfortunately, computer calculations have shown that
$\Stwo(9,4)$ and $\Stwo(10,5)$ are \textbf{not} lattices,
rendering this lattice-theoretic method of proof invalid for $d \geq 4$
(and resolving negatively Conjecture 2.13 of~\cite{EdelmanReiner1996}).
However we would still conjecture the following:

\begin{Conjecture}
The image of $\phi:\oo \rightarrow \intbar\Stwo(n,d)$ is
exactly the subposet consisting of those closed intervals in $\Stwo(n,d)$ 
whose open interval is not contractible.
\end{Conjecture}

\noindent
As in Section~\ref{sec:baues}, this conjecture would resolve
in the affirmative the Baues problem for triangulations of all
cyclic polytopes. It is easy to see that one direction in this conjecture is
true, namely that any interval in the image of $\phi$ is isomorphic
to a Cartesian product of posets isomorphic to $\Stwo(n_\alpha,d)$
for various $n_\alpha$, and hence has proper part homotopy equivalent
to a sphere.  Consequently, the above conjecture also has as
a corollary the calculation of the homotopy type and
M\"obius function for all (open) intervals in $\Stwo(n,d)$.

\item Do the partial orders $\Sone(n,d), \Stwo(n,d)$ coincide?

\end{enumerate}
 
\section{Acknowledgments}
The authors thank P. Webb for pointing out Lemma~\ref{Webb-lemma},
which gives an easier and more conceptual proof than their original proof of 
Lemma~\ref{lemma:Baues:atomic-proper-intervals}.

%%%%%%%%%%%%%%%%%%%%%%%%%%%%%%%%%%%%%%%%%%%%%%%%%%%%%%%%%%%%%%%%%%%%%%%%%%%%%%%
%% References:                                                               %%
%% 05/09/96 jr                                                               %%
%%%%%%%%%%%%%%%%%%%%%%%%%%%%%%%%%%%%%%%%%%%%%%%%%%%%%%%%%%%%%%%%%%%%%%%%%%%%%%%
\providecommand{\bysame}{\leavevmode\hbox to3em{\hrulefill}\thinspace}

\end{document}

%% file: bauclic_zigzagexample.pstex_t
\begin{picture}(0,0)%
\special{psfile=bauclic_zigzagexample.pstex}%
\end{picture}%
\setlength{\unitlength}{0.00087500in}%
\begingroup\makeatletter\ifx\SetFigFont\undefined
% extract first six characters in \fmtname
\def\x#1#2#3#4#5#6#7\relax{\def\x{#1#2#3#4#5#6}}%
\expandafter\x\fmtname xxxxxx\relax \def\y{splain}%
\ifx\x\y   % LaTeX or SliTeX?
\gdef\SetFigFont#1#2#3{%
  \ifnum #1<17\tiny\else \ifnum #1<20\small\else
  \ifnum #1<24\normalsize\else \ifnum #1<29\large\else
  \ifnum #1<34\Large\else \ifnum #1<41\LARGE\else
     \huge\fi\fi\fi\fi\fi\fi
  \csname #3\endcsname}%
\else
\gdef\SetFigFont#1#2#3{\begingroup
  \count@#1\relax \ifnum 25<\count@\count@25\fi
  \def\x{\endgroup\@setsize\SetFigFont{#2pt}}%
  \expandafter\x
    \csname \romannumeral\the\count@ pt\expandafter\endcsname
    \csname @\romannumeral\the\count@ pt\endcsname
  \csname #3\endcsname}%
\fi
\fi\endgroup
\begin{picture}(4860,958)(114,-322)
\put(3826,-286){\makebox(0,0)[b]{\smash{\SetFigFont{12}{14.4}{rm}(b)}}}
\put(1351,-286){\makebox(0,0)[b]{\smash{\SetFigFont{12}{14.4}{rm}(a)}}}
\put(136,254){\makebox(0,0)[lb]{\smash{\SetFigFont{12}{14.4}{rm}$S_1$}}}
\put(136, 29){\makebox(0,0)[lb]{\smash{\SetFigFont{12}{14.4}{rm}$S_2$}}}
\put(497,478){\makebox(0,0)[b]{\smash{\SetFigFont{12}{14.4}{rm}1}}}
\put(722,478){\makebox(0,0)[b]{\smash{\SetFigFont{12}{14.4}{rm}2}}}
\put(947,478){\makebox(0,0)[b]{\smash{\SetFigFont{12}{14.4}{rm}3}}}
\put(1172,478){\makebox(0,0)[b]{\smash{\SetFigFont{12}{14.4}{rm}4}}}
\put(1397,478){\makebox(0,0)[b]{\smash{\SetFigFont{12}{14.4}{rm}5}}}
\put(1622,478){\makebox(0,0)[b]{\smash{\SetFigFont{12}{14.4}{rm}6}}}
\put(1847,478){\makebox(0,0)[b]{\smash{\SetFigFont{12}{14.4}{rm}7}}}
\put(2072,478){\makebox(0,0)[b]{\smash{\SetFigFont{12}{14.4}{rm}8}}}
\put(2297,478){\makebox(0,0)[b]{\smash{\SetFigFont{12}{14.4}{rm}9}}}
\put(2973,477){\makebox(0,0)[b]{\smash{\SetFigFont{12}{14.4}{rm}1}}}
\put(3198,477){\makebox(0,0)[b]{\smash{\SetFigFont{12}{14.4}{rm}2}}}
\put(3423,477){\makebox(0,0)[b]{\smash{\SetFigFont{12}{14.4}{rm}3}}}
\put(3648,477){\makebox(0,0)[b]{\smash{\SetFigFont{12}{14.4}{rm}4}}}
\put(3873,477){\makebox(0,0)[b]{\smash{\SetFigFont{12}{14.4}{rm}5}}}
\put(4098,477){\makebox(0,0)[b]{\smash{\SetFigFont{12}{14.4}{rm}6}}}
\put(4323,477){\makebox(0,0)[b]{\smash{\SetFigFont{12}{14.4}{rm}7}}}
\put(4548,477){\makebox(0,0)[b]{\smash{\SetFigFont{12}{14.4}{rm}8}}}
\put(4773,477){\makebox(0,0)[b]{\smash{\SetFigFont{12}{14.4}{rm}9}}}
\put(2611, 29){\makebox(0,0)[lb]{\smash{\SetFigFont{12}{14.4}{rm}$S_2$}}}
\put(2611,254){\makebox(0,0)[lb]{\smash{\SetFigFont{12}{14.4}{rm}$S_1$}}}
\put(496,254){\makebox(0,0)[b]{\smash{\SetFigFont{12}{14.4}{rm}$*$}}}
\put(946,254){\makebox(0,0)[b]{\smash{\SetFigFont{12}{14.4}{rm}$*$}}}
\put(721, 29){\makebox(0,0)[b]{\smash{\SetFigFont{12}{14.4}{rm}$*$}}}
\put(1171, 29){\makebox(0,0)[b]{\smash{\SetFigFont{12}{14.4}{rm}$*$}}}
\put(1621,254){\makebox(0,0)[b]{\smash{\SetFigFont{12}{14.4}{rm}$*$}}}
\put(2071,254){\makebox(0,0)[b]{\smash{\SetFigFont{12}{14.4}{rm}$*$}}}
\put(2296,254){\makebox(0,0)[b]{\smash{\SetFigFont{12}{14.4}{rm}$*$}}}
\put(2071, 29){\makebox(0,0)[b]{\smash{\SetFigFont{12}{14.4}{rm}$*$}}}
\put(1846, 29){\makebox(0,0)[b]{\smash{\SetFigFont{12}{14.4}{rm}$*$}}}
\put(1621, 29){\makebox(0,0)[b]{\smash{\SetFigFont{12}{14.4}{rm}$*$}}}
\put(2971,254){\makebox(0,0)[b]{\smash{\SetFigFont{12}{14.4}{rm}$*$}}}
\put(3421,254){\makebox(0,0)[b]{\smash{\SetFigFont{12}{14.4}{rm}$*$}}}
\put(4096,254){\makebox(0,0)[b]{\smash{\SetFigFont{12}{14.4}{rm}$*$}}}
\put(4546,254){\makebox(0,0)[b]{\smash{\SetFigFont{12}{14.4}{rm}$*$}}}
\put(4771,254){\makebox(0,0)[b]{\smash{\SetFigFont{12}{14.4}{rm}$*$}}}
\put(4546, 29){\makebox(0,0)[b]{\smash{\SetFigFont{12}{14.4}{rm}$*$}}}
\put(4321, 29){\makebox(0,0)[b]{\smash{\SetFigFont{12}{14.4}{rm}$*$}}}
\put(4096, 29){\makebox(0,0)[b]{\smash{\SetFigFont{12}{14.4}{rm}$*$}}}
\put(3421, 29){\makebox(0,0)[b]{\smash{\SetFigFont{12}{14.4}{rm}$*$}}}
\put(3196, 29){\makebox(0,0)[b]{\smash{\SetFigFont{12}{14.4}{rm}$*$}}}
\end{picture}

%% file: bauclic_gapexample.pstex_t
\begin{picture}(0,0)%
\special{psfile=bauclic_gapexample.pstex}%
\end{picture}%
\setlength{\unitlength}{0.00087500in}%
\begingroup\makeatletter\ifx\SetFigFont\undefined
% extract first six characters in \fmtname
\def\x#1#2#3#4#5#6#7\relax{\def\x{#1#2#3#4#5#6}}%
\expandafter\x\fmtname xxxxxx\relax \def\y{splain}%
\ifx\x\y   % LaTeX or SliTeX?
\gdef\SetFigFont#1#2#3{%
  \ifnum #1<17\tiny\else \ifnum #1<20\small\else
  \ifnum #1<24\normalsize\else \ifnum #1<29\large\else
  \ifnum #1<34\Large\else \ifnum #1<41\LARGE\else
     \huge\fi\fi\fi\fi\fi\fi
  \csname #3\endcsname}%
\else
\gdef\SetFigFont#1#2#3{\begingroup
  \count@#1\relax \ifnum 25<\count@\count@25\fi
  \def\x{\endgroup\@setsize\SetFigFont{#2pt}}%
  \expandafter\x
    \csname \romannumeral\the\count@ pt\expandafter\endcsname
    \csname @\romannumeral\the\count@ pt\endcsname
  \csname #3\endcsname}%
\fi
\fi\endgroup
\begin{picture}(4860,958)(114,-322)
\put(3826,-286){\makebox(0,0)[b]{\smash{\SetFigFont{12}{14.4}{rm}(b)}}}
\put(1351,-286){\makebox(0,0)[b]{\smash{\SetFigFont{12}{14.4}{rm}(a)}}}
\put(497,478){\makebox(0,0)[b]{\smash{\SetFigFont{12}{14.4}{rm}1}}}
\put(722,478){\makebox(0,0)[b]{\smash{\SetFigFont{12}{14.4}{rm}2}}}
\put(947,478){\makebox(0,0)[b]{\smash{\SetFigFont{12}{14.4}{rm}3}}}
\put(1172,478){\makebox(0,0)[b]{\smash{\SetFigFont{12}{14.4}{rm}4}}}
\put(1397,478){\makebox(0,0)[b]{\smash{\SetFigFont{12}{14.4}{rm}5}}}
\put(1622,478){\makebox(0,0)[b]{\smash{\SetFigFont{12}{14.4}{rm}6}}}
\put(1847,478){\makebox(0,0)[b]{\smash{\SetFigFont{12}{14.4}{rm}7}}}
\put(2072,478){\makebox(0,0)[b]{\smash{\SetFigFont{12}{14.4}{rm}8}}}
\put(2297,478){\makebox(0,0)[b]{\smash{\SetFigFont{12}{14.4}{rm}9}}}
\put(2973,477){\makebox(0,0)[b]{\smash{\SetFigFont{12}{14.4}{rm}1}}}
\put(3198,477){\makebox(0,0)[b]{\smash{\SetFigFont{12}{14.4}{rm}2}}}
\put(3423,477){\makebox(0,0)[b]{\smash{\SetFigFont{12}{14.4}{rm}3}}}
\put(3648,477){\makebox(0,0)[b]{\smash{\SetFigFont{12}{14.4}{rm}4}}}
\put(3873,477){\makebox(0,0)[b]{\smash{\SetFigFont{12}{14.4}{rm}5}}}
\put(4098,477){\makebox(0,0)[b]{\smash{\SetFigFont{12}{14.4}{rm}6}}}
\put(4323,477){\makebox(0,0)[b]{\smash{\SetFigFont{12}{14.4}{rm}7}}}
\put(4548,477){\makebox(0,0)[b]{\smash{\SetFigFont{12}{14.4}{rm}8}}}
\put(4773,477){\makebox(0,0)[b]{\smash{\SetFigFont{12}{14.4}{rm}9}}}
\put(497, 28){\makebox(0,0)[b]{\smash{\SetFigFont{12}{14.4}{rm}$o$}}}
\put(721,254){\makebox(0,0)[b]{\smash{\SetFigFont{12}{14.4}{rm}$e$}}}
\put(722, 28){\makebox(0,0)[b]{\smash{\SetFigFont{12}{14.4}{rm}$o$}}}
\put(947,253){\makebox(0,0)[b]{\smash{\SetFigFont{12}{14.4}{rm}$e$}}}
\put(1172,253){\makebox(0,0)[b]{\smash{\SetFigFont{12}{14.4}{rm}$e$}}}
\put(1397, 28){\makebox(0,0)[b]{\smash{\SetFigFont{12}{14.4}{rm}$o$}}}
\put(1622, 28){\makebox(0,0)[b]{\smash{\SetFigFont{12}{14.4}{rm}$o$}}}
\put(1847,253){\makebox(0,0)[b]{\smash{\SetFigFont{12}{14.4}{rm}$e$}}}
\put(3197,253){\makebox(0,0)[b]{\smash{\SetFigFont{12}{14.4}{rm}$e$}}}
\put(4322,253){\makebox(0,0)[b]{\smash{\SetFigFont{12}{14.4}{rm}$e$}}}
\put(4322, 28){\makebox(0,0)[b]{\smash{\SetFigFont{12}{14.4}{rm}$o$}}}
\put(4547, 28){\makebox(0,0)[b]{\smash{\SetFigFont{12}{14.4}{rm}$o$}}}
\put(496,254){\makebox(0,0)[b]{\smash{\SetFigFont{12}{14.4}{rm}$*$}}}
\put(946, 29){\makebox(0,0)[b]{\smash{\SetFigFont{12}{14.4}{rm}$*$}}}
\put(1171, 29){\makebox(0,0)[b]{\smash{\SetFigFont{12}{14.4}{rm}$*$}}}
\put(1396,254){\makebox(0,0)[b]{\smash{\SetFigFont{12}{14.4}{rm}$*$}}}
\put(1621,254){\makebox(0,0)[b]{\smash{\SetFigFont{12}{14.4}{rm}$*$}}}
\put(1846, 29){\makebox(0,0)[b]{\smash{\SetFigFont{12}{14.4}{rm}$*$}}}
\put(2071, 29){\makebox(0,0)[b]{\smash{\SetFigFont{12}{14.4}{rm}$*$}}}
\put(2071,254){\makebox(0,0)[b]{\smash{\SetFigFont{12}{14.4}{rm}$*$}}}
\put(2296,254){\makebox(0,0)[b]{\smash{\SetFigFont{12}{14.4}{rm}$*$}}}
\put(3196, 29){\makebox(0,0)[b]{\smash{\SetFigFont{12}{14.4}{rm}$*$}}}
\put(3421, 29){\makebox(0,0)[b]{\smash{\SetFigFont{12}{14.4}{rm}$*$}}}
\put(3421,254){\makebox(0,0)[b]{\smash{\SetFigFont{12}{14.4}{rm}$*$}}}
\put(4096,254){\makebox(0,0)[b]{\smash{\SetFigFont{12}{14.4}{rm}$*$}}}
\put(4501,254){\makebox(0,0)[b]{\smash{\SetFigFont{12}{14.4}{rm}$*$}}}
\put(4771,254){\makebox(0,0)[b]{\smash{\SetFigFont{12}{14.4}{rm}$*$}}}
\put(4771, 29){\makebox(0,0)[b]{\smash{\SetFigFont{12}{14.4}{rm}$*$}}}
\put(4096, 29){\makebox(0,0)[b]{\smash{\SetFigFont{12}{14.4}{rm}$*$}}}
\put(2296, 29){\makebox(0,0)[b]{\smash{\SetFigFont{12}{14.4}{rm}$*$}}}
\put(136,254){\makebox(0,0)[lb]{\smash{\SetFigFont{12}{14.4}{rm}$S_1$}}}
\put(136, 29){\makebox(0,0)[lb]{\smash{\SetFigFont{12}{14.4}{rm}$S_2$}}}
\put(2611,254){\makebox(0,0)[lb]{\smash{\SetFigFont{12}{14.4}{rm}$S_1$}}}
\put(2611, 29){\makebox(0,0)[lb]{\smash{\SetFigFont{12}{14.4}{rm}$S_2$}}}
\end{picture}

%% file: bauclic_slide.pstex_t
\begin{picture}(0,0)%
\epsfig{file=bauclic_slide.pstex}%
\end{picture}%
\setlength{\unitlength}{0.00087500in}%
\begingroup\makeatletter\ifx\SetFigFont\undefined
% extract first six characters in \fmtname
\def\x#1#2#3#4#5#6#7\relax{\def\x{#1#2#3#4#5#6}}%
\expandafter\x\fmtname xxxxxx\relax \def\y{splain}%
\ifx\x\y   % LaTeX or SliTeX?
\gdef\SetFigFont#1#2#3{%
  \ifnum #1<17\tiny\else \ifnum #1<20\small\else
  \ifnum #1<24\normalsize\else \ifnum #1<29\large\else
  \ifnum #1<34\Large\else \ifnum #1<41\LARGE\else
     \huge\fi\fi\fi\fi\fi\fi
  \csname #3\endcsname}%
\else
\gdef\SetFigFont#1#2#3{\begingroup
  \count@#1\relax \ifnum 25<\count@\count@25\fi
  \def\x{\endgroup\@setsize\SetFigFont{#2pt}}%
  \expandafter\x
    \csname \romannumeral\the\count@ pt\expandafter\endcsname
    \csname @\romannumeral\the\count@ pt\endcsname
  \csname #3\endcsname}%
\fi
\fi\endgroup
\begin{picture}(4319,1416)(172,-748)
\put(3634,-253){\makebox(0,0)[lb]{\smash{\SetFigFont{12}{14.4}{rm}\special{ps: gsave 0 0 0 setrgbcolor}$T$\special{ps: grestore}}}}
\put(3634,-709){\makebox(0,0)[lb]{\smash{\SetFigFont{12}{14.4}{rm}\special{ps: gsave 0 0 0 setrgbcolor}$i(f(T))$\special{ps: grestore}}}}
\put(213,-709){\makebox(0,0)[lb]{\smash{\SetFigFont{12}{14.4}{rm}\special{ps: gsave 0 0 0 setrgbcolor}$T$\special{ps: grestore}}}}
\put(2494,-709){\makebox(0,0)[lb]{\smash{\SetFigFont{12}{14.4}{rm}\special{ps: gsave 0 0 0 setrgbcolor}$f(T)$\special{ps: grestore}}}}
\end{picture}